\documentclass[12pt,a4paper]{article}

\usepackage[latin1]{inputenc}
\usepackage{color}

\usepackage{bussproofs}
\usepackage{graphicx}
\usepackage{graphics}
\usepackage{latexsym}
\usepackage{amssymb,amsmath}
\usepackage{euscript}
\usepackage{enumerate}
\usepackage{amsmath}
\usepackage{amsfonts}
\usepackage{multicol}
\usepackage{etex}
\usepackage[all]{xy}
\usepackage{color}
\usepackage{graphicx}
\usepackage{tikz}
\usepackage[latin1]{inputenc}

\graphicspath{ {images/} }

\usepackage[colorlinks=true, pdfstartview=FitV, linkcolor=blue, citecolor=blue, urlcolor=blue, pagebackref,]{hyperref}

\newtheorem{teo}{Theorem}[section]
\newtheorem{coro}[teo]{Corollary}
\newtheorem{lema}[teo]{Lemma}
\newtheorem{prop}[teo]{Proposition}

\newtheorem{defi}[teo]{Definition}

\newtheorem{rem}[teo]{Remark}
\newtheorem{thm}{Theorem}[section]

\newtheorem{lem}[thm]{Lemma}

\newtheorem{defn}[thm]{Definition}

\newenvironment{dem}{\noindent {\em Proof:} }{\hfill $\square$ \bigskip}
\newenvironment{prueba}[1]{\noindent\bf Proof of #1. \rm}{$\quad \hfill \square$ \bigskip}
\numberwithin{equation}{section}

\newcommand{\R}{\mathbb{R}}

\newcommand{\N}{\mathbb{N}}
\renewcommand{\H}{\mathbb{H}}

\newcommand{\supp}{\mathop{\text{\rm supp}}}
\newcommand{\alphat}{\tilde{\alpha}}
\newcommand{\deltat}{\tilde{\delta}}
\newcommand{\menort}{\lesssim}
\newcommand{\mayort}{\gtrsim}
\newcommand{\zR}{\mathbb{R}}
\newcommand{\zN}{\mathbb{N}}

\newcommand{\abs}[1]{\left|#1\right|}

\newcommand{\norm}[1]{\left\|#1\right\|}

\def\XXint#1#2#3{{\setbox0=\hbox{$#1{#2#3}{\int}$}
    \vcenter{\hbox{$#2#3$}}\kern-.5\wd0}}

\usepackage{fancyhdr}

\fancypagestyle{encabezadotesis}{
\fancyhead{}
\fancyhead[RO]{\em \nouppercase{\rightmark}}
\fancyhead[LE]{\bf \nouppercase{\leftmark}}
\cfoot{\thepage}
}

\begin{document}
\title{On two weighted problems for commutators of classical operators with optimal behaviour on the parameters involved and extrapolation results}
\author{Gladis Pradolini, Wilfredo Ramos and Jorgelina Recchi}

%
%
%
%
%
%

\date{}
\maketitle

\renewcommand{\thefootnote}{\fnsymbol{footnote}}
\footnotetext{
The  authors are supported by CONICET and
UNL, UNNE and UNS respectively.
}

\begin{abstract}
We give two-weighted norm estimates for higher order commutator of classical operators such as singular integral and fractional type operators, between weighted $L^p$ and certain spaces that include Lipschitz, BMO and Morrey spaces. We also give the optimal parameters involved with these results, where the optimality is understood in the sense that the parameters defining the corresponding spaces belong to certain region out of which the classes of weights are satisfied by trivial weights. We also exhibit pairs of non-trivial weights in the optimal region satisfying the conditions required. Finally, we exhibit an extrapolation result that allows us to obtain boundedness
results of the type described above in the variable setting and for a great variety of operators, by starting from analogous inequalities in the classical context. In order to get this result we prove a Calder\'on-Scott type inequality with weights that connects adequately the spaces involved.
\end{abstract}

\section{Introduction}

There is a vast evidence of the effect that the continuity properties of several operators from Harmonic Analysis have on the regularity estimations of the solution of certain partial differential equation, where the latter benefit from the former. Among a great variety of operators appear the commutators of singular and fractional integral operators, which have been extensively studied. We can cite, as examples, \cite{BCe}, \cite{ST2}, \cite{CUF}, \cite{DLZ}, \cite{LK}, \cite{P}, \cite{P2}, \cite{PP}, \cite{PR2}, \cite{HST}, \cite{Chanillo} and \cite{Janson} between a great amount of papers on this topic. For example, in \cite{BCe} the authors extend a result from \cite{CRW}, by proving the boundedness of commutators of parabolic Calderón-Zygmund operators between $L^p$ spaces. The main motivation arose when they were proving $L^p$ estimates for parabolic equations with VMO coefficients. They also use this result in \cite{BCM} for some related problems associated to ultraparabolic operators with VMO coefficients, (see also \cite{BraCe}, \cite{BraCeru}, \cite{ChFL}, \cite{ChFL2} and \cite{Rios} for related problems).

More specifically, in \cite{DPR} the authors obtained one-weight estimates for integral operators of fractional type and their higher
order commutators between certain spaces. The results include weighted $L^p$-Lipschitz$(\beta)$ estimates, where $\beta\ge 0$ depends on some parameter related with the symbol of the commutator. The relation between $p$ and $\beta$ is standard, that is, $\beta/n=\theta/n-1/p$, where $0<\theta<n$. The Lipschitz spaces considered in that article are the same as those given in \cite{Pra}, which are generalizations of some weighted version of the $BMO$ spaces introduced by Muckenhoupt and Wheeden in \cite{MW}.
In \cite{MP} similar problems were studied for singular integral type operators and their higher order commutators with Lipschitz symbols and standard relations between the parameters involved.

In this paper we consider the optimal parameters that appear in connection with two-weighted problems of the type described above, but the range of $p$ and $\beta$ goes over the region where the classes of weights involved are satisfied by pairs $(v,w)$ of non-trivial weights. The Lipschitz spaces are the same as in \cite{DPR} and \cite{MP}, but we can also prove the results for BMO and Morrey spaces, extending some previous theorems proved in \cite{Pra} for the fractional integral operator. Concretely, we prove two-weighted norm estimates for singular integral and fractional type operators and their higher order commutators between weighted $L^p$ and certain generalized spaces related to a parameter $\beta$. We also give the optimal parameters involved with these results, where the optimality is understood in the sense that the parameters $p$ and $\beta$ belong to certain region out of which the classes of weights are satisfied by trivial weights, that is $w=0$ or $v=\infty$. Moreover, we exhibit concrete pairs of non-trivial weights in the optimal region satisfying the conditions required on the weights, where the boundedness results includes values of $\beta$ describing Lipschitz, BMO and Morrey spaces, that is, $0<\beta<1$, $\beta=0$ and $\beta<0$, respectively. Our results extend those contained in \cite{Pra} for the fractional integral operator. We also prove that a one-weight result can only holds whenever the relation between the parameters is the standard.

Related with the extrapolation results, in 1982, J. L.  Rubio de Francia (\cite{R}) proved that the
$\mathcal{A}_{p}$ classes enjoy a very interesting extrapolation
property: If for some $1\leq p_{0}<\infty$, an operator preserves
$L^{p_{0}}(w)$ for any $w\in \mathcal{A}_{p_{0}}$, then necessarily
preserves the  $L^{p}(w)$ space for every $1<p<\infty$ and every
$w\in\mathcal{A}_{p}$. This result was very surprising and
established the begining for a great variety of research articles in
this topic. For example, in 1988, E. Harboure, R.  Mac\'ias and  C.
Segovia (\cite{HMS}) proved that the $\mathcal{A}(p,q)$ classes have
a similar extrapolation property, that is, if $T$ is a sublinear
operator such that the inequality
\begin{equation}\label{a10}
\norm{Tf\,w}_{q}\leq C\norm{fw}_{p}
\end{equation}
holds for some pair $(p_{0},q_{0})$,  $1<p_{0}\leq q_{0}<\infty$ and
every $w\in\mathcal{A}(p_{0},q_{0})$, then  \eqref{a10} holds for
every pair $(p,q)$, $1<p\leq q<\infty$, satisfying the condition
$$1/p-1/q=1/p_{0}-1/q_{0}$$ and every $w\in \mathcal{A}(p,q)$.
Moreover, they also proved that this property is not only exclusive
for the boundedness between weighted Lebesgue spaces, but also it is
possible to extrapolate based on a continuity behavior from weighted
Lebesgue spaces into weighted BMO spaces. In other words, if
$1<\beta<\infty$ and $T$ is a sublinear operator satisfying
\begin{equation}\label{a11}
\norm{Tf}_{\mathbb{L}_w (0)}\leq C\norm{f/w}_{\beta},
\end{equation}
for every ball $B\subset \mathbb{R}^{n}$ and every weight $w^{-1}\in \mathcal{A}(\beta,\infty)$; then,  if $1<p<\beta$, $1/p-1/q=1/\beta$ and $w^{-1}\in\mathcal{A}(p,q)$, the inequality \eqref{a10}, with $w$ replaced by $w^{-1}$, holds for the pair  $(p,q)$  provided that the left hand side is finite.
On the other hand, in \cite{CPR} the authors give extrapolation results that allow them
to obtain continuity properties of certain operators, of the type $L^p_w- L^q_w$ and $L^p_w-\mathbb L_w(\tilde{\delta})$, starting with hypothesis of continuity of the type $L^s_w-\mathbb L_w({\delta})$ for some related parameters.

In the scale of variable Lebesgue spaces, some extrapolation results
for classical operators were given in \cite{CUFMP} and in \cite{CPR}
in order to obtain unweighted results in the spirit of the previous
ones described above. Later, in \cite{CUW} the authors give
boundedness results between weighted variable Lebesgue spaces
starting from similar inequalities but in the classical context.

In this paper we also prove an extrapolation result that allows us to obtain boundedness results for a great variety of operators in Harmonic Analysis between weighted variable Lebesgue spaces into certain weighted variable Lipschitz spaces starting from analogous inequalities in the classical context, that is, with constant parameters. In order to achieve this result we also prove a Calder\'on-Scott type inequality with weights that generalizes the inequality given in (\cite{CSt}) and that allows us to connect variable Lebesgue spaces with variable Lipschitz spaces.

The paper is organized as follows. In \S $\ref{preliminaries}$ we give the preliminaries in order to state the first part of our main results. The study of the properties of the class of weights is given in \S $\ref{propweight}$, where we also give the optimal range of the parameters involved. Later, in \S $\ref{sprimeraparte}$ we give the proofs of the results established in \S $\ref{preliminaries}$. Finally, in \S $\ref{extra1}$, \S $\ref{extra2}$ and \S $\ref{extra3}$ we state the extrapolations results mentioned in the introduction, some applications and the corresponding proofs, respectively.

\section{Preliminaries and main results}\label{preliminaries}

In this section we give the definitions of the operators we shall be dealing with and the functional class of the symbols in order to define the commutators.

Let $0<\delta< 1$. We say that a function $b$ belongs to the space $\Lambda(\delta)$ if there exists a positive constant $C$ such that, for every $x,y \in \R^n$
$$
|b(x)-b(y)|\leq  C |x-y|^\delta.
$$
The smallest of such constants will be denoted by $\|b\|_{\Lambda(\delta)}$. The space $\Lambda(\delta)$ is the well known Lipschitz space in the classical literature. We shall be dealing with commutators  with symbols belonging to this class of functions.

We say that a weight $u$ belongs to the reverse H\"older class $RH(s)$, $1\le s<\infty$ if there exists a positive constant $C$ such that
\begin{equation}\label{reverse}
\left(\frac{1}{|B|}\int_B u^s(x)dx\right)^{1/s}\leq C\, \frac{u(B)}{|B|}.
\end{equation}
If $s=\infty$ we say that $u\in RH(\infty)$ if (\ref{reverse}) holds with the average of the left hand side replaced by $\|\chi_B u\|_{\infty}$.


We say that $\Phi:[0,\infty)\to [0,\infty)$ is a Young function if it is increasing, convex and verifies $\Phi(0)=0$ and $\Phi(t)\to \infty$ when $t\to \infty$. The $\Phi$-Luxemburg average of a locally integrable function $f$ over a ball $B$ is defined by
$$
\|f\|_{\Phi,B}= \inf \left\{\lambda>0: \frac{1}{|B|}\int_B \Phi\left(\frac{|f(x)|}{\lambda}\right)\, dx \leq 1 \right\}.
$$

Given two functions $f$ and $g$, by $\menort $ and $\mayort$ we shall mean that there exists a positive constant $c$ such that $f\leq c\,g$ and $cf\ge g$, respectively. We say that $f \simeq g$ if both inequalities $f\menort g$ and $f\mayort g$ hold.

If $\Phi$ is a Young function, the following H\"older's type inequality holds for every pair of measurable functions $f$, $g$
\begin{equation}\label{Holder}
\frac{1}{|B|}\int_B |f(x)g(x)|dx\menort \|f\|_{\Phi,B}\|g\|_{\tilde{\Phi},B},
\end{equation}
where $\tilde{\Phi}$ is the complementary Young function of $\Phi$, defined by
$$
\tilde{\Phi}(t)=\sup_{s>0}\{st-\Phi(s)\}.
$$
It is not difficult to see that $t\leq\Phi^{-1}(t)\tilde{\Phi}^{-1}(t)\leq 2t$ for every $t>0$. Moreover, given $\Phi$, $\Psi$ and $\Theta$ Young functions verifying that $\Phi^{-1}(t)\Psi^{-1}(t)\menort \Theta^{-1}(t)$ for every $t>0$, the following generalization of (\ref{Holder}) holds
$$
\|fg\|_{\Theta,B}\menort\|f\|_{\Phi,B}\|g\|_{\Psi,B}.
$$
For more information about this topics, see \cite{RR}.

We shall consider singular integral type operators $T$ which are the convolution with a kernel $K$, that is $T$ is bounded on $L^2(\R^n)$ and if $x\notin \supp f$
\begin{equation} \label{integral singular}
Tf(x)=\int_{\R^n}K(x-y)f(y)dy.
\end{equation}
The kernel $K$ is a measurable function defined away from $0$.

We shall also consider  fractional operators of convolution type $T_\alpha$, $0< \alpha<n$, defined by

\begin{equation}\label{operador fraccionario}
T_\alpha f(x)=\int_{\R^n}K_\alpha(x-y)f(y)dy,
\end{equation}
where the kernel $K_\alpha$ is not identically zero.
We supposse that both kernels, $K$ and $K_{\alpha}$ in (\ref{integral singular}) and (\ref{operador fraccionario}) respectively, satisfy certain size and smoothness conditions that we shall define as follows.

We say that the  kernel $K$ satisfies the size condition $S_{0}^*$ if there exists a positive constant $C$ such that
$$
|K(x)|\leq \frac{C}{\, \,|x|^{n}},
$$
and we say that $K$ satisfies the smoothness condition $H^{*}_{0, \infty}$ if there exists a positive constant $C$ and $0<\gamma\leq 1$ such that
$$
|K(x-y)-K(x'-y)|+|K(y-x)-K(y-x')|\le C \frac{|x-x'|^\gamma}{\, \, \, \, |x-y|^{n+\gamma}},
$$
whenever $|x-y|\geq 2|x-x'|$.

Let $0<\alpha<n$. We say that the  kernel $K_\alpha$ satisfies the size condition $S_{\alpha}^*$ if there exists a positive constant $C$ such that
$$
|K_\alpha(x)|\leq \frac{C}{\, \,|x|^{n-\alpha}},
$$
and we say that $K_\alpha$ satisfies the smoothness condition $H^{*}_{\alpha, \infty}$ if there exists a positive constant $C$ and $0<\gamma\leq 1$ such that
$$
|K_\alpha(x-y)-K_\alpha(x'-y)|+|K_\alpha(y-x)-K_\alpha(y-x')|\le C \frac{|x-x'|^\gamma}{\, \, \, \, |x-y|^{n-\alpha+\gamma}},
$$
whenever $|x-y|\geq 2|x-x'|$.

It is easy to check that the Hilbert transform $H$, with kernel $K(x)=1/x$, and the fractional integral operator $I_\alpha$, with kernel $K_\alpha(x)=|x|^{\alpha-n}$, satisfy conditions $S^{*}_{\alpha}$  and  $H^{*}_{\alpha, \infty}$ for $\alpha=0$ and $0<\alpha<n$, respectively.

Let $\mathcal{T}$ denotes both, the singular integral type operator $T$ or the fractional type operator $T_{\alpha}$. We can formally define the commutator with symbol $b\in L^1_{\mathrm {loc}}(\R^n)$, by
$$
[b,\mathcal{T}]f=b\,\mathcal{T}f-\mathcal{T}(bf).
$$
The commutator of order $m\in \N \cup \{0\}$ of $\mathcal{T}$ is defined by
$$
\mathcal{T}_b^0=\mathcal{T}, \;\;\;\;\; \mathcal{T}^m_b=[b,\mathcal{T}^{m-1}_b].
$$
We denote $\mathcal{T}^m_b={T}^m_b$ and $\mathcal{T}^m_b=T_{\alpha,b}^m$ the commutators with symbol $b$ of the singular integral operator $T$ and the fractional type operator $T_{\alpha}$, respectively.

As we have said, we are  interested in studying the boundedness properties of the commutators $T_{b}^m$  and $T_{\alpha,b}^m$ with symbol $b\in \Lambda(\delta)$,  from weighted Lebesgue spaces into weighted versions of certain spaces, including Lipschitz, BMO and Morrey type spaces. For $\beta \in \R$ and a weight $w$, these spaces are denoted by $\mathbb{L}_w(\beta)$ and collect the functions $f\in L^{1}_{\mathrm {loc}}(\R^n)$ that satisfy
$$
\|f\|_{\mathbb{L}_w(\beta)}=\sup_{B}\frac{\|(1/w)\chi_B\|_{\infty}}{|B|^{1+\beta}}\int_B |f(x)-m_B(f)|\,dx <\infty.
$$

When $\beta=0$, $\mathbb{L}_w(0)$ is a weighted version of the bounded  mean oscillation  space introduced by Muckenhoupt and Wheeden in \cite{MW}. Moreover,  $\mathbb{L}_1(\beta)$ gives the known Lipschitz integral space for $\beta$ in the range $0<\beta<1/n$ and the Morrey space, for $-1<\beta<0$.

\smallskip

Troughout this paper we denote $\alphat=m\delta +\alpha$, where $0\leq\alpha<n$, $m\in \N \cup \{0\}$ and  $0<\delta <1$.

We now introduce a class of pairs of weights. As we shall see, they are related with the boundedness results described above.

\begin{defi}\label{definicion de la clase de pesos}
Let  $0\leq\alpha<n$, $0<\delta<\min \{\gamma,(n-\alpha)/m\}$ and $1< r\leq \infty$. Let $\deltat\leq\delta$. We say that a pair of weights $(v,w)$ belongs to $\H (r,\alphat,\deltat)$, if
\begin{equation}\label{condicion del peso}
\sup_B\|(1/w)\chi_B\|_{\infty}|B|^{\frac{\delta-\deltat}{n}}\left(\int_{\R^n}\frac{v^{r'}(y)}{\left(|B|^{1/n}+|x_B-y|\right)^{r'(n-\alphat+\delta)}}\,dy\right)^{1/{r'}} <\infty,
\end{equation}
where the supremun is taken over every ball $B$ with center $x_B$.
In the case $r=1$, $(\ref{condicion del peso})$ should be understood as
\begin{equation}\label{condicion del pesoinf}
\sup_B\|(1/w)\chi_B\|_{\infty}|B|^{\frac{\delta-\deltat}{n}}\left\|\frac{v(\cdot)}{\left(|B|^{1/n}+|x_B-\cdot|\right)^{(n-\alphat+\delta)}}\right\|_{\infty} <\infty,
\end{equation}
\end{defi}

When $w=v$ we say $w\in \H(r,\alphat,\deltat)$ if

\begin{equation}\label{condicion para un solo peso}
\sup_B\|(1/w)\chi_B\|_{\infty}|B|^{\frac{\delta-\deltat}{n}}\left(\int_{\R^n}\frac{w^{r'}(y)}{\left(|B|^{1/n}+|x_B-y|\right)^{r'(n-\alphat+\delta)}}\,dy\right)^{1/{r'}} <\infty,
\end{equation}
with the obvious changes when $r=1$.

When $\delta=1$ and $m=0$, the class above was introduced in \cite{Pra}. In that paper the author proves that the corresponding condition characterizes the boundedness of the fractional integral operator of order $\gamma$, $I_{\gamma}$, between $L^{p}(v)$ and $\mathbb{L}_{w}(\deltat/n)$ spaces. Moreover, when $v=w$ some similar problems were studied in \cite{HSV2}.

The following lemma shows, as we expect, that if $w\in \H(r,\alphat,\deltat)$, then the parameters involved satisfy a relation that is a generalization of the well known one that appears in the case $m=0$ (see, for example, \cite{Pra}).

\begin{lema} \label{un solo peso solo en el borde}
Let  $0\leq\alpha<n$, $0<\delta<\min \{\gamma,(n-\alpha)/m\}$ and $1\le r\leq \infty$. Let $\deltat\leq \delta$. If $w\in \H (r,\alphat,\deltat)$, then $\deltat=\alphat-n/r$.

\end{lema}
\begin{dem}
Let $1<r\le \infty$. Since $w\in \H (r,\alphat,\deltat)$, using (\ref{condicion para un solo peso}) we have,
\begin{align*}
\inf_B w&\gtrsim |B|^{\frac{\delta-\deltat}{n}}\left(\int_{B}\frac{w^{r'}(y)}{\left(|B|^{1/n}+|x_B-y|\right)^{r'(n-\alphat+\delta)}}\,dy\right)^{1/{r'}}\\
& \mayort |B|^{\frac{\delta-\deltat}{n}-(\frac{1}{r}-\frac{\alphat}{n}+\frac{\delta}{n})}\left(\frac{1}{|B|}\int_{B}w^{r'}(y)\,dy\right)^{1/{r'}}\\
& \mayort |B|^{-\frac{\deltat}{n}-\frac{1}{r}+\frac{\alphat}{n}}\inf_B w.
\end{align*}
Then, this inequality holds if $\deltat=\alphat-n/r$. The case $r=1$ follows with the corresponding changes.
\end{dem}

\begin{rem}In view of the lemma above,  when we consider one weight, we rename the classes $\H (r,\alphat,\deltat)$ as $\H (r,\alphat)$ in order to show their dependence on the parameter $r$ and $\alphat$.
\end{rem}
\smallskip

We are now in a position to state our main results. We first state them for singular integral operators with the corresponding pairs if weights belonging to $\H(r,m\delta,\deltat)$, that is $\alpha=0$ in Definition \ref{definicion de la clase de pesos}. Recall that $K$ satisfies conditions $S^*_0$ and $H_{0,\infty}^{*}$.
\color{black}
\begin{teo} \label{teo para integrales singulares}
Let $0<\delta<\min \{\gamma,n/m\}$, $1\leq r \le \infty$ and $\deltat\leq m\delta-n/r$. Let $b\in \Lambda(\delta)$. If $(v,w) \in \H(r,m\delta,\deltat)$, then the inequality
$$
\|T^m_b f\|_{\mathbb{L}_w(\deltat/n)} \menort\,\|b\|^m_{\Lambda(\delta)}\left\| f/v \right\|_{L^r(\R^n)},
$$
holds for every $f$ such that $f/v \in L^r(\R^n)$.
\end{teo}
As a consequence of Theorem \ref{teo para integrales singulares} and Lemma \ref{un solo peso solo en el borde}, applied to the case $\alpha=0$, we obtain the following corollary, which was proved in \cite{MP}. In fact, we shall see that the classes $\H(r,m\delta)$ are the same as the classes $A_{r,\infty}$ considered in that paper.

\begin{coro}\label{luchigla} Let $0<\delta<\min \{\gamma,n/m\}$, $1\leq r \le \infty$ and $\deltat= m\delta-n/r$. Let $b\in \Lambda(\delta)$. If $w \in \H(r,m\delta)$, then the inequality
$$
\|T^m_b f\|_{\mathbb{L}_w(\deltat/n)} \menort\,\|b\|^m_{\Lambda(\delta)}\left\| f/w \right\|_{L^r(\R^n)},
$$
holds for every $f$ such that $f/w \in L^r(\R^n)$.
\end{coro}
When $m=1$ and $w=1$ the corollary above was proved in \cite{MY} in the general setting of non-doubling measures. On the other hand, when $m=0$ this result is a generalization of that obtained in \cite{MW2} for the Hilbert transform.

\smallskip

We now state the main results for fractional integral operators. Recall that $K_\alpha$ satisfies condition $S^*_\alpha$ and $H_{\alpha,\infty}^{*}$.

\begin{teo}\label{teorema A}
Let $0<\alpha<n$, $0<\delta<\min\{\gamma, (n-\alpha)/m\}$, $1\leq r \le \infty$ and $\deltat\leq m\delta-n/r$. Let $b\in \Lambda(\delta)$. If $(v,w)\in \H(r,\alphat,\deltat)$, then the inequality
$$
\|T^m_{\alpha,b}f\|_{\mathbb{L}_w(\deltat/n)} \menort\,\|b\|^m_{\Lambda(\delta)}\left\| f/v \right\|_{L^r(\R^n)},
$$
holds for every $f$ such that $f/v \in L^r(\R^n)$.
\end{teo}

When $T_{\alpha}$ is the fractional integral operator of order $\alpha$, $I_\alpha$, the result above was proved in \cite{Pra} for the case $m=0$.

As a consequence of Theorem \ref{teorema A} and Lemma \ref{un solo peso solo en el borde} we have the following result, which was previously proved in \cite{DPR}.

\begin{coro}\label{corolario fraccionario}
Let $0<\alpha<n$, $0<\delta<\min\{\gamma, (n-\alpha)/m\}$, $1\leq r \le \infty$ and $\deltat=\alphat-n/r$. Let $b\in \Lambda(\delta)$. If $w\in \H(r,\alphat)$ then the inequality
$$
\|T^m_{\alpha,b}f\|_{\mathbb{L}_w(\deltat/n)} \menort\,\|b\|^m_{\Lambda(\delta)}\left\| f/w \right\|_{L^r(\R^n)},
$$
holds for every $f$ such that $f/w \in L^r(\R^n)$.
\end{coro}

When $r=n/\alpha$ and  $\deltat=0$, the corollary above was obtained in \cite{MW}. On the other hand, a similar result was proved in \cite{MY} for $m=1$, $w=1$ and $\deltat=\delta+\alpha-n/r$ in the general context of non-doubling measures.

\smallskip


%
%
%
%
%
%

Related with the classes $\H(r,\alphat,\deltat)$ we have obtained an extrapolation theorem that allows us to obtain boundedness results for a great variety of operators in Harmonic Analysis between weighted variable Lebesgue spaces into certain weighted variable Lipschitz spaces, starting from analogous inequalities in the classical context, that is, with constant parameters. In order to state this result, we give some preliminaries definitions and notations.

Let $p:\zR^n\rightarrow[1,\infty)$ be a mesurable function. For $A\subset\zR^n$ we define
	$p^{-}_{A}=\text{ess}\inf_{x\in A}p(x)$ {and} $
	p^{+}_{A}=\text{ess}\sup_{x\in A}p(x).$
	For simplicity we denote $p^{+}=p^{+}_{\zR^n}$ and
	$p^{-}=p^{-}_{\zR^n}$.
	
	We say that $p\in\mathcal{P}(\zR^n)$ if  $1<p^{-}\le p(x)\le
	p^{+}<\infty$ for almost every $x\in\zR^n$.
	Let $p\in\mathcal{P}(\zR^n)$. The variable exponent Lebesgue space $L^{p(\cdot)}(\zR^{n})$ is the set of the measurable functions $f$ defined on $\zR^n$ such that, for some positive $\lambda$, the  convex functional modular
	\begin{equation*}
	\varrho(f/\lambda)=\int_{\zR^n}|f(x)/\lambda|^{p(x)}\,dx
	\end{equation*}
	is finite. A Luxemburg norm can be defined in $L^{p(\cdot)}(\zR^{n})$ by taking
	\begin{equation*}
	 \norm{f}_{L^{p(\cdot)}(\zR^{n})}=\inf\{\lambda>0:\varrho(f/\lambda)\le1\}.
	\end{equation*}
	
	The variable Lebesgue spaces are special cases of
	Museliak-Orlicz spaces (see \cite{M}), and generalize the classical
	Lebesgue spaces. For more information see, for example, \cite{CUF2},
	\cite{DHHR} or \cite{KR}.
	
	We say that $p\in\mathcal{P}^{\text{log}}(\zR^n)$
	if  $p\in\mathcal{P}(\zR^n)$ and  it satisfies the following
	inequalities
	\begin{equation}
	|p(x)-p(y)|\le \frac{C}{\log(e+1/|x-y|)},\hspace{1cm}\text{for
		every}~x,y\in\zR^n,
	\end{equation}
	and
	\begin{equation}
	|p(x)-p(y)|\le \frac{C}{\log(e+|x|)},\hspace{1cm}\text{with}~|y|\ge |x|.
	\end{equation}
	
	Let $w$ be a weight, $1<\beta<\infty$ and
	$p\in\mathcal{P}(\mathbb{R}^{n})$ such that $\beta\leq p^{-}\leq
	p(x)\leq p^{+}<\frac{n\beta}{(n-\beta)_{+}}$, where $\rho_+$ is defined as $\rho$ if $\rho>0$ and $0$ is $\rho\leq 0$. Let
	${\delta(x)}/{n}={1}/{\beta}-{1}/{p(x)}$. The space
	$\mathbb{L}_{w}(\delta(\cdot))$ is defined by the set of the measurable
	functions $f$ such that
	\begin{equation*}
	 |||f|||_{\mathbb{L}_{w}(\delta(\cdot))}=\sup_{B}\frac{||(1/w)\chi_{B}||_{\infty}}{\abs{B}^{{1}/{\beta}}\norm{\chi_{B}}_{p'(\cdot)}}\int_{B}\abs{f-m_{B}f}<\infty.
	\end{equation*}
	
	When $p(x)$ is equal to a constant $p$, this space coincides with the
	space $\mathbb{L}_{w}({n}/{\beta}-{n}/{p})$ defined above. When $w=1$, the space $\mathbb{L}_1(\delta(\cdot))$ was defined in
	\cite{RSV}.
	
	\medskip

	Let $p(\cdot),~
	q(\cdot)\in\mathcal{P}(\zR^{n})$. We say that a
	weight $u\in A_{p(\cdot),q(\cdot)}$ if there exists a positive
	constant $C$ such that for every ball $B$
	\begin{equation*}\|u\chi_B\|_{q(\cdot)}\|u^{-1}\chi_B\|_{p'(\cdot)}\le C\|\chi_{B}\|_{q(\cdot)}\|\chi_{B}\|_{p'(\cdot)}.
	\end{equation*}
	When $p(\cdot)\in \mathcal{P}^{\log}(\zR^n)$ and $p(\cdot)=q(\cdot)$
	this is the $A_{p(\cdot)}$ class introduced in \cite{CUDH}. It is
	well known that, if $p(\cdot)\in \mathcal{P}^{\log}(\zR^n)$ then
	$u\in A_{p(\cdot)}$ if and only if $M: L_u^{p(\cdot)}
	\hookrightarrow L_u^{p(\cdot)}$ (\cite{CUDH}). On the other hand, when $p(\cdot)=p$ and $q(\cdot)=q$ are constant exponents, the $A_{p,q}$ class was introduced in \cite{MW}.
	
	We also say that  $u\in A_{p(\cdot),\infty}$ if there exists a positive constant $C$ such that
	\begin{equation*}
	\|u\chi_{B}\|_{\infty}\|u^{-1}\chi_{B}\|_{p'(\cdot)}\le C\|\chi_{B}\|_{p'(\cdot)}.
	\end{equation*}

When $q=\infty$, and $p(\cdot)=p$ is a constant exponent, we say that ${u}\in A_{p,\infty}$ if ${u}^{-p'}
\in A_1$, (see \cite{MW} for more information about these classes).

	\medskip
	
	We are now in a position to state our main theorem related with
	extrapolation results from extreme values of the exponents involved.

	\begin{teo}\label{teo2}
		Let $1<\beta<\infty$, $0\leq \delta <\min\{{\color{blue}n/\beta}, 1\}$ and let $s>1$ be
		such that ${1}/{s}={1}/{\beta}-{\delta}/{n}$. Let $p(\cdot)\in
		\mathcal{P}^{\log}(\zR^n)$ such that
		${1}/{p(\cdot)}={1}/{\beta}-{\tilde{\delta}(\cdot)}/{n}$, with
		$0\leq \tilde{\delta}(\cdot)\leq \delta$ and $0<\nu<p^-$. If $f$
		and $g$ are two measurable functions such that the inequality
		$$|||f|||_{\mathbb{L}_{w}(\delta)}\leq C \|gw\|_{s} $$
		holds for every weight $w$ such that $w^{-\nu}\in A_{{s}/{\nu}, \infty}$, then there exists a positive constant $C$ such that the inequality
		$$|||f|||_{\mathbb{L}_{w}(\tilde{\delta}(\cdot))} \leq C \|gw\|_{p(\cdot)} $$
		holds for every $w$ such that $w^{-\nu}\in A_{{p(\cdot)}/{\nu},\, \infty}$.
	\end{teo}
	
	\begin{rem}
		When $w=1$ and $\nu=1$ the result above was proved in {\rm \cite{CPR}}.
	\end{rem}
	
	
	Let $1<\beta<\infty$ and  $s(\cdot)\in\mathcal{P}(\zR^{n})$ be  such
	that ${\delta(x)}/{n}={1}/{\beta}-{1}/{s(x)}\geq 0$, the
	$\delta(\cdot)$$\,$-$\,$sharp maximal operator $f^{\sharp}_{\delta(\cdot)}$ is defined by
	$$f^{\sharp}_{\delta(\cdot)}(x)=\sup_{B\ni x}\frac{1}{\abs{B}^{1/\beta}\norm{\chi_{B}}_{s'(\cdot)}}\int_{B}\abs{f-m_{B}f},$$
	where the supremum is taking over every ball $B$ contaning $x$.
	
	\smallskip
	
	In order to achieved Theorem \ref{teo2} we prove a Calder\'on-Scott type inequality (\cite{CSt}) with weights that allows us to connect variable Lebesgue spaces with variable Lipschitz spaces. An
	unweighted version can be found in \cite{CPR}.
	\begin{prop}\label{P4}
		Let $0\leq \tilde{\delta}(\cdot)< \delta(\cdot)<1$. Let $s(\cdot)$
		and $p(\cdot)\in \mathcal{P}^{\log}(\zR^n)$ be defined by
		$\delta(\cdot)/n=1/\beta-1/s(\cdot)$ and
		$\tilde{\delta}(\cdot)/n=1/\beta-1/p(\cdot)$, respectively. If $w^{-1}\in
		A_{{n}/{(\delta(\cdot) -\tilde{\delta}(\cdot))},\infty}$ then
		$$\|f^{\sharp}_{\tilde{\delta}(\cdot)}/w\|_{\infty}\leq C\|f^{\sharp}_{\delta(\cdot)}/w\|_{{n}/{(\delta -\tilde{\delta}(\cdot)})}.$$
	\end{prop}



%
%

\section{Properties of the classes of weights}\label{propweight}

We  give some properties of the  classes of weights $\H (r,\alphat,\deltat)$. Recall that $\deltat\leq \alphat-n/r$ and $\alphat=m\delta+\alpha$, where $0\leq \alpha< n$.

\begin{lema}\label{lema del doble v} Let $0<\alpha<n$, $0<\delta<\min\{\gamma, (n-\alpha)/m\}$, $1\leq r \le \infty$ and $\deltat\leq \alphat-n/r$.
If $(v,w) \in\H (r,\alphat,\deltat)$ then
$$
\|v\chi_{2B}\|_{r'}\menort |B|^{1+\frac{(\deltat-\alphat)}{n}}\inf_Bw,
$$
where the constant involved is independent of $B$  and $w$.
\end{lema}
\begin{dem}
Since $(v,w)\in \H(r, \alphat,\deltat)$, we know that
$$
|B|^{(\delta-\deltat)/n}\left(\int_{\R^n}\frac{v^{r'}(y)}{(|B|^{1/n}+|x_B-y|)^{r'(n-\alphat+\delta)}}dy\right)^{1/{r'}}\leq C\inf_B
$$
for every ball $B\subset\R^n$, where $x_B$ is the center of $B$. Then, we have
\begin{align*}
v^{r'}(2B)& =\frac{|B|^{(n-\alphat+\delta)r'/n}}{|B|^{(n-\alphat+\delta)r'/n}}\int_{2B}v^{r'}(y)dy\\
& \menort |B|^{(n-\alphat+\delta)r'/n}\int_{2B}\frac{v^{r'}(y)}{(|B|^{1/n}+|x_B-y|)^{r'(n-\alphat+\delta)}}dy\\
& \menort |B|^{(n-\alphat+\delta)r'/n}\int_{\R^n}\frac{v^{r'}(y)}{(|B|^{1/n}+|x_B-y|)^{r'(n-\alphat+\delta)}}dy\\
& \menort |B|^{(n-\alphat+\delta)r'/n}\left(\frac{\inf_Bw}{|B|^{\frac{\delta-\deltat}{n}}}\right)^{r'}\\
& \menort |B|^{(1+(\deltat-\alphat)/n)r'}\left(\inf_Bw\right)^{r'}.
\end{align*}
\end{dem}

%
%

\begin{lema}\label{lema condicion de peso es equivalente a la parte local mas la global}
Let  $0\leq \alpha<n$, $0<\delta<\min\{\gamma, (n-\alpha)/m\}$ and  $1\le r\leq \infty$. The condition $\H (r,\alphat,\deltat)$ is equivalent to  the  following two inequalities
\begin{equation}\label{parte local}
\|(1/w)\chi_B\|_{\infty}|B|^{\frac{\alphat-\deltat}{n}-\frac{1}{r}}\left(\frac{1}{|B|}\int_B v^{r'}(y)dy\right)^{1/{r'}}\menort \, 1
\end{equation}
and
\begin{equation}\label{parte global}
\displaystyle\|(1/w)\chi_B\|_{\infty}|B|^{\frac{\delta-\deltat}{n}}\left(\int_{\R^n-B}\frac{v^{r'}(y)}{|x_B-y|^{r'(n-\alphat+\delta)}}dy\right)^{1/{r'}}\menort 1
\end{equation}
hold simultaneously for every ball $B\subset \R^n$, where $x_B$ is the center of $B$.
\end{lema}

\begin{dem} Let us first see that $\H (r,\alphat,\deltat)$ implies (\ref{parte local}) and (\ref{parte global}). Given $B=B(x_B,R)$, since for $y\in B$ we have $|x_B-y|\menort |B|^{1/n}$, from condition $\H (r,\alphat,\deltat)$ we obtain that
$$
|B|^{\frac{\alphat-\deltat}{n}-\frac{1}{r}}\left(\frac{1}{|B|}\int_B v^{r'}(y)dy\right)^{1/{r'}}\menort \, \inf_Bw
$$
which proves (\ref{parte local}). In order to prove (\ref{parte global}), let us observe that if $y\in \R^n\setminus B$ then $|x_B-y|\mayort |B|^{1/n}$. Consequently, we have
$$
\displaystyle|B|^{\frac{\delta-\deltat}{n}}\left(\int_{\R^n-B}\frac{v^{r'}(y)}{|x_B-y|^{r'(n-\alphat+\delta)}}dy\right)^{1/{r'}}\menort \inf_Bw.
$$
Assume (\ref{parte local}) and (\ref{parte global}), then for every ball $B$
\begin{align*}
\|(1/w)\chi_B\|_{\infty}|B|^{\frac{\delta-\deltat}{n}}&\left(\int_{\R^n}\frac{v^{r'}(y)}{\left(|B|^{1/n}+|x_B-y|\right)^{r'(n-\alphat+\delta)}}\,dy\right)^{1/{r'}} \\
&\menort \|(1/w)\chi_B\|_{\infty}|B|^{\frac{\alphat-\deltat}{n}-\frac{1}{r}}\left(\frac{1}{|B|}\int_{B}v^{r'}\right)^{1/{r'}} \\
&\quad + \|(1/w)\chi_B\|_{\infty}|B|^{\frac{\delta-\deltat}{n}}\left(\int_{\R^n\setminus B}\frac{v^{r'}(y)}{\left(|B|^{1/n}+|x_B-y|\right)^{r'(n-\alphat+\delta)}}\,dy\right)^{1/{r'}}\\
&\menort 1.
\end{align*}
\end{dem}

\begin{lema} \label{local implica global}
Let $0\le \alpha<n$, $0<\delta<\min\{\gamma, (n-\alpha)/m\}$, $1 \leq r\le \infty$ and  $\deltat< \delta$. Then the local condition $(\ref{parte local})$ implies the global condition $(\ref{parte global})$.
\end{lema}

\begin{dem}
Let $B=B(x_B,R)$ and $B_k=2^kB$. By using the local condition \ref{parte local}, we have that
\begin{align*}
&\frac{\|(1/w)\chi_B\|_{\infty}}{\displaystyle|B|^{\frac{\deltat-\delta}{n}}}\left(\int_{\R^n-B}\frac{v^{r'}(y)}{|x_B-y|^{r'(n-\alphat+\delta)}}dy\right)^{1/{r'}}\\
& \hspace{3cm} \menort \frac{\|(1/w)\chi_B\|_{\infty}}{\displaystyle|B|^{\frac{\deltat-\delta}{n}}}\sum_{k=1}^{\infty}|B_k|^{-1+\frac{\alphat-\delta}{n}}\left(\int_{B_k-B_{k-1}} v^{r'}(y)dy\right)^{1/{r'}}\\
&\hspace{3cm}\le \frac{1}{\displaystyle|B|^{\frac{\deltat-\alphat}{n}+\frac{1}{r}}}\sum_{k=1}^{\infty}2^{k({{\alphat-\delta}-\frac{n}{r}})}\|(1/w)\chi_{B_k}\|_{\infty}\left(\frac{1}{|B_k|}\int_{B_k} v^{r'}(y)dy\right)^{1/{r'}}\\
&\hspace{3cm}\menort  \sum_{k=1}^{\infty}2^{k(\deltat-\delta)}\le C.
\end{align*}

\end{dem}

%
%

When $\deltat=\delta$ the condition $\H(r,\alphat,\deltat)$ can not be reduced to (\ref{parte local}).  In order to prove it, we first recall the following well known estimates whose proof is left to the reader.

\begin{lema}\label{lema 3.5 de Gladis}
Let $B=B(x_B,R)\subset \R^n$ and $\alpha >-n$. The following statements holds.
\begin{enumerate}
\item If $|x_B|\leq R$,  $\int_B |x|^{\alpha}dx\approx R^{\alpha+n}$
\item If $|x_B|>R$,  $\int_B|x|^{\alpha}dx\approx|x_B|^{\alpha}R^n$.
\end{enumerate}
\end{lema}
 We now prove the assert above. Let $\alphat-n<\deltat=\delta<\alphat-n/r$. Let $(v,w)$ be the pair of weights defined by
$$
w\equiv1 \;\;\;\;\;\;\;\;\;\; \text{and}\;\;\;\;\;\; \;\;\;\; v(x)=|x|^{n/r-\alphat+\delta}.
$$
We shall see that $(v,w)$ satisfies (\ref{parte local}) but does not satisfy (\ref{parte global}).

Given $B=B(x_B, R)$, $R>0$, from Lemma \ref{lema 3.5 de Gladis}, it's easy to check that
$$
|B|^{\frac{\alphat-\deltat}{n}-\frac{1}{r}}\left(\frac{1}{|B|}\int_B v^{r'}(y)dy\right)^{1/{r'}}\menort \, \inf_Bw.
$$
We shall now see that $(v,w)$ does not satisfy (\ref{parte global}). Taking $B=B(0,R)$ we have
\begin{align*}
\left(\int_{\R^n\setminus B}\frac{v^{r'}(y)}{|x_B-y|^{r'(n-\alphat+\delta)}}dy\right)^{1/{r'}}&= \left(\int_{|y|>R}\frac{|y|^{r'(n/r-\alphat+\delta)}}{|y|^{r'(n-\alphat+\delta)}}dy\right)^{1/{r'}}\\
& =\left(\int_{|y|>R}\frac{1}{|y|^n}dy\right)^{1/{r'}},
\end{align*}
and this last integral is infinite.\\
It is important to emphasize that there exists examples of non-trivial pairs of weights  that satisfy, for the case $\deltat=\delta$, both conditions (\ref{parte local}) and (\ref{parte global}), that  is, pairs that belong to $\H(r,\alphat,\delta)$.

We shall now study the ranges of $r$ and $\deltat$ for which the pairs of weights that verify $\H(r,\alphat,\deltat)$ are trivial, that is $v=0$ a.e. or $w=\infty$ a.e. Recall that $\alphat=m \delta + \alpha$.

{\color{black}
\begin{prop}
Let $0\le\alpha<n$, $0<\delta<\min\{\gamma, (n-\alpha)/m\}$ and $1\leq r\le \infty$.
\begin{enumerate}
\item\label{i} If $\deltat >\delta$ or $\deltat>\alphat-n/r$, then $(v,w)\in \H(r,\alphat,\deltat)$ if and only if $v=0$ almost everywhere in $\R^n$.
\item\label{ii}   If $\deltat=\alphat-n/r=\delta$ the same conclusion as in $i)$ holds.
\item\label{iii} If $\deltat<\alphat-n$ then $(v,w)\in \H(r,\alphat,\deltat)$ if and only if $v=0$ a.e. in $\R^n$ or $w=\infty$ a.e. in $\R^n$.
\end{enumerate}
\end{prop}
\begin{dem}
Let us first prove $(\ref{i})$ when $\deltat >\delta$. Let $B=B(x,R)$ where $x$ is a Lebesgue point of $w$. Then, since $(v,w)\in \H(r,\alphat,\deltat)$
$$
\left(\int_{\R^n}\frac{v^{r'}(y)}{(|B|^{1/n}+|x-y|)^{r'(n-\alphat+\delta)}}dy\right)^{1/{r'}}\menort \inf_Bw \, |B|^{(\deltat-\delta)/n}\menort \frac{w(B)}{|B|} |B|^{(\deltat-\delta)/n}.
$$
Letting $R \to 0$, the right hand side of the inequality above tends to zero. Then, we obtain that
\begin{equation*}
\left(\int_{\R^n}\frac{v^{r'}(y)}{|x-y|^{r'(n-\alphat+\delta)}}dy\right)^{1/{r'}}=0
\end{equation*}
and, consequently, $v=0$ almost everywhere.\\
Let us now consider $\deltat>\alphat-n/r$. Since $(v,w)\in \H(r,\alphat,\deltat)$, by Lemma \ref{lema condicion de peso es equivalente a la parte local mas la global} we have that
$$
\left(\frac{v^{r'}(B)}{|B|}\right)^{1/r'}\menort \inf_Bw \, |B|^{1/r-(\alphat-\deltat)/n}\menort \frac{w(B)}{|B|}|B|^{1/r-(\alphat-\deltat)/n},
$$
where we choose $B(x,R)$ as in the case above. Thus we get
$$
\lim_{R\to 0} \frac{w(B)}{|B|}|B|^{1/r-(\alphat-\deltat)/n}=0,
$$
and so
$$
\limsup_{R\to 0} \frac{v^{r'}(B(x,R))}{|B(x,R)|}=0.
$$
Then, by standar arguments, we can deduce that $v(x)=0$ in almost everywhere $x\in \R^n$.

If $r=\infty$  we must consider $1/r=0$ and $r'=1$ in the previous proofs.

In order to estimate $(\ref{ii})$, suppose that $\deltat=\delta=\alphat-n/r$. If $\alphat >\delta$ then $r=n/{(\alphat-\delta)}$ and, if $\alphat=\delta$ then $r=\infty$. Let $(v,w)\in \H(r,\alphat,\deltat)$ and $=B(x_0, R)\subset \R^{n}$, where $x_0$ will be choosen later. Since $\deltat=\delta$ we get that
$$
\left(\int_{\R^n}\frac{v^{r'}(y)}{(|B|^{1/n}+|x_0-y|)^{r'(n-\alphat+\delta)}}dy\right)^{1/{r'}}\menort \inf_Bw \le \frac{w(B)}{|B|}.
$$
Observe that $r'< \infty$ since the case $r=1$ is not reached. Since  $n-\alphat+\delta=n/r'$, the inequality above can be rewritten as

\begin{equation}\label{A}
\left(\int_{\R^n}\frac{v^{r'}(y)}{(|B|^{1/n}+|x_0-y|)^{n}}dy\right)^{1/{r'}}\menort\frac{w(B)}{|B|}.
\end{equation}
Let us now assume that the set $E=\{x: v(x)>0\}$ have positive measure. Since
$$
E=\cup_{k}E_k=\cup_k \{x\in E: v(x)>1/k\}
$$
then there exists $k_0$ such that $|E_{k_0}|>0$. We select $x_0\in E_{k_0}$ a Lebesgue point of $w$ and a density point of $E_{k_0}$ such that $w(x_0)$ is finite.
Thus, from (\ref{A}) and letting $R \to 0$, we get
$$
\left(\int_{\R^n}\frac{v^{r'}(y)}{|x_0-y|^{n}}dy\right)^{1/{r'}}\menort w(x_0),
$$
and the same conclusion holds for almost everywhere $x_0 \in E_{k_0}$. Then,
\begin{equation}\label{3.18}
\left(\int_{\R^n}\frac{v^{r'}(y)}{|x_0-y|^{n}}dy\right)^{1/{r'}} < \infty
\end{equation}
for almost everywhere $x_0 \in E_{k_0}$. On the other hand we have
\begin{equation}\label{3.19}
\left(\int_{\R^n}\frac{v^{r'}(y)}{|x_0-y|^{n}}dy\right)^{1/{r'}}>\frac{1}{k_0}\left(\int_{B\cap E_{k_0}}\frac{1}{|x_0-y|^{n}}dy\right)^{1/{r'}}.
\end{equation}
By following standar arguments it is not difficult to prove that the right hand side tends to $\infty$ when  $R \to 0$,  which contradicts (\ref{3.18}). Therefore $v=0$ almost everywhere.

We now proceed with the proof of $(\ref{iii})$. Let $B=B(x_B,R)$ and let $B_0$ be a ball such that $B_0\subset B$. From condition $(\ref{parte local})$ we have that
\begin{align*}
\|(1/w)\chi_{B_0}\|_{\infty}\left(\int_{B_0}v^{r'}(y)dy\right)^{1/{r'}}\le \|(1/w)\chi_{B}\|_{\infty}\left(\int_{B}v^{r'}(y)dy\right)^{1/{r'}}\menort \, R^{\deltat - \alphat+n}.
\end{align*}
Since $\deltat<\alphat-n$ the right hand side tends to $0$ when $R\to \infty$ and so, the statement holds.
\end{dem}}

The next result exhibits the range of  $r$ and $\deltat$ for which the class  $\H(r,\alphat,\deltat)$ is satisfied by non-trivial pair of weights. See also Figure $\ref{rango para clase H}$ for a better understanding.
\includegraphics{grafico}
\begin{teo} \label{teo donde viven los pesos}
Let $0\le \alpha<n$ and $0<\delta<\min\{\gamma, (n-\alpha)/m\}$. Then there exist non-trivial pairs of weights $(v,w)$ that verify the condition $\H(r,\alphat,\deltat)$ in the range of $r$ and $\deltat$ given by
$$
1\leq r\leq \infty \;\;\;\;\;\;\;\text{and}\;\;\;\;\;\;\;\alphat-n\le \deltat \leq \min\{\delta, \alphat-n/r\}
$$
excluding the case $\deltat=\delta$ when $\alphat-n/r=\delta$.
\end{teo}

\begin{dem}
	We shall divide the proof in the following cases
	\begin{enumerate}[a.]
		\item\label{primero} $\deltat=\delta<\alphat-n/r$,
		\item\label{segundo} $\alphat-n<\deltat<\delta\le\alphat-n/r$
		\item\label{tercero} $\alphat-n<\deltat\le\alphat-n/r<\delta$
		\item\label{cuarto} $\deltat=\alphat-n$, $1\le r\le\infty$.	
	\end{enumerate}

Let us first prove $(\ref{primero})$. Let $\theta$ and $\beta$ be such that
\begin{equation*}
\alphat-n/r-\delta<\theta<n/r',\quad {\rm and} \quad 0<\beta=\theta+\delta+n/r-\alphat.
\end{equation*}
Let $(v,w)$ be the pair of weights defined by
\begin{equation*}
w(x)=|x|^{-\beta}  \qquad v(x)=|x|^{-\theta}.
\end{equation*}
Let us see that $(v,w)\in \H(r, \alphat,\deltat)$, which is equivalent to prove that there exists a positive constant $C$ such that both inequalities
\begin{equation}\label{aa}
\|(1/w)\chi_B\|_{\infty}|B|^{\frac{\alphat-\delta}{n}-\frac{1}{r}}\left(\frac{1}{|B|}\int_B v^{r'}(y)dy\right)^{1/{r'}}\le C
\end{equation}
and
\begin{equation}\label{bb}
\|(1/w)\chi_B\|_{\infty}\left(\int_{\R^n-B}\frac{v^{r'}(y)}{|x_B-y|^{r'(n-\alphat+\delta)}}dy\right)^{1/{r'}}\le C
\end{equation}
hold.
Let us first estimate $(\ref{aa})$. If $|x_B|\le R$, then $\|(1/w)\chi_B\|_{\infty}\menort R^{\beta}$ and, by lemma $\ref{lema 3.5 de Gladis}$ and the fact that  $\beta<\theta$, we obtain that
\begin{equation*}
\|(1/w)\chi_B\|_{\infty}|B|^{\frac{\alphat-\delta}{n}-\frac{1}{r}}
\left(\frac{1}{|B|}\int_B v^{r'}(y)dy\right)^{1/{r'}} \menort R^{\beta+\alphat-\delta-n/r-\theta}\menort 1.
\end{equation*}
On the other hand, if $|x_B|> R$ then $\|(1/w)\chi_B\|_{\infty}\menort |x_B|^{\beta}$. Thus, from Lemma  $\ref{lema 3.5 de Gladis}$ we get
\begin{align*}
\|(1/w)\chi_B\|_{\infty}|B|^{\frac{\alphat-\delta}{n}-\frac{1}{r}}
\left(\frac{1}{|B|}\int_B v^{r'}(y)dy\right)^{1/{r'}} &\menort {|x_B|^{\beta-\theta}}R^{\alphat-\delta-n/r}\\
&\menort R^{\beta-\theta+\theta-\beta}\menort 1.
\end{align*}

We now proceed with the estimation of $(\ref{bb})$. Let $B_k=2^kB$, $k\in \N$. Then
\begin{align}
I&=\|(1/w)\chi_B\|_{\infty}\left(\int_{\R^n-B}\frac{v^{r'}(y)}{|x_B-y|^{r'(n-\alphat+\delta)}}dy\right)^{1/{r'}}\nonumber\\&\label{ddd}\le \|(1/w)\chi_B\|_{\infty}\sum_{k=1}^{\infty}\frac{1}{(2^{k}R)^{n-\alphat+\delta}}\left(\int_{B_k}v^{r'}\right)^{1/r'}
\end{align}
Let us first assume $|x_B|\le R$. By Lemma $\ref{lema 3.5 de Gladis}$ we get
\begin{equation*}
I \menort R^{\beta-n+\alphat-\delta-\theta+n/r'}\sum_{k=1}^{\infty}2^{-k\beta}\menort 1.
\end{equation*}
If $|x_B|> R$, there exists $N_1\in \N$ such that $2^{N_1}R< |x_B|\le2^{N_1+1}R$. Then, $(\ref{ddd})$ can be written as $S_1+S_2$, where
\begin{equation*}
S_1=\|(1/w)\chi_B\|_{\infty}\sum_{k=1}^{N_1}\frac{1}{(2^{k}R)^{n-\alphat+\delta}}\left(\int_{B_k}v^{r'}\right)^{1/r'}
\end{equation*}
and
\begin{equation*}S_2=\|(1/w)\chi_B\|_{\infty}\sum_{k=N_1+1}^{\infty}\frac{1}{(2^{k}R)^{n-\alphat+\delta}}\left(\int_{B_k}v^{r'}\right)^{1/r'}.
\end{equation*}
We first estimate $S_1$. Since $|x_B|>2^kR$ for $k=1,2,...,N_1$, then, from Lemma $\ref{lema 3.5 de Gladis}$ we obtain that
\begin{align*}
S_1&\menort |x_B|^{\beta} R^{\alphat-n-\delta}\sum_{k=1}^{N_1}2^{-k(n-\alphat+\delta)}|x_B|^{-\theta}(2^{k}R)^{n/r'}\\
&\menort R^{\theta-\beta}|x_B|^{\beta-\theta}\sum_{k=1}^{N_1}2^{k(\theta-\beta)}\\
&\menort R^{\theta-\beta}|x_B|^{\beta-\theta} |x_B|^{\theta-\beta}  R^{\beta-\theta} \menort 1.
\end{align*}
In order to estimate $S_2$ recall that if $k>N_1$ then $|x_B|\le 2^kR$. Thus, by Lemma  $\ref{lema 3.5 de Gladis}$ we obtain that
\begin{align*}
S_2&\menort |x_B|^{\beta}R^{\alphat-n-\delta}\sum_{k=N_1+1}^{\infty}2^{-k(n-\alphat+\delta)}(2^{k}R)^{-\theta+n/r'}\\
&\menort \left(\frac{|x_B|}{R}\right)^{\beta}\sum_{k=N_1+1}^{\infty}2^{-k\beta}
\cong \left(\frac{|x_B|}{R}\right)^{\beta} \left(\frac{R}{|x_B|}\right)^{\beta}\menort 1.
\end{align*}
The proof for the case $\deltat=\delta<\alphat-n/r$ is concluded.

In the cases  $(\ref{segundo})$, $(\ref{tercero})$ and $(\ref{cuarto})$ we have that $\deltat<\delta$ and so, by Lemma $\ref{local implica global}$ we only have to estimate condition $(\ref{parte local})$.

Let us begin with the case $(\ref{segundo})$, that is $\alphat-n<\deltat<\delta\le\alphat-n/r$. Let $\theta$ and $\beta$ be such that
\begin{equation*}
\alphat-n/r-\deltat<\theta<n/r',\quad {\rm and} \quad 0<\beta=\theta+\deltat+n/r-\alphat.
\end{equation*}
Let $(v,w)$ be the pair of weights defined by
\begin{equation*}
w(x)=|x|^{-\beta}  \qquad v(x)=|x|^{-\theta}.
\end{equation*}
Then $(v,w)\in \H(r, \alphat,\deltat)$. From Lemma  $\ref{lema 3.5 de Gladis}$, if $|x_B|\le R$ we get
\begin{equation*}
\|(1/w)\chi_B\|_{\infty}|B|^{\frac{\alphat-\deltat}{n}-\frac{1}{r}}\left(\frac{1}{|B|}\int_B v^{r'}(y)dy\right)^{1/{r'}}\menort R^{\beta+\alphat-\deltat-n/r-\theta}\menort 1,
\end{equation*}
and, if $|x_B|> R$
\begin{align*}
\|(1/w)\chi_B\|_{\infty}|B|^{\frac{\alphat-\deltat}{n}-\frac{1}{r}}\left(\frac{1}{|B|}\int_B v^{r'}(y)dy\right)^{1/{r'}}&\menort |x_B|^{\beta}R^{\alphat-\deltat-n/r}|x_B|^{-\theta}\menort\left(\frac{R}{|x_B|}\right)^{\theta-\beta}\menort 1,
\end{align*}
since $\theta-\beta>0$.

Let us now consider the situation in $(\ref{tercero})$, that is $\alphat-n<\deltat\le\alphat-n/r<\delta$. Let $(v,w)$ be the pair of weights defined by
$$w(x)=1 \quad {\rm and }\quad v(x)=|x|^{-\theta}, $$
where $\theta=\alphat-n/r-\deltat$. Then $(v,w)\in \H(r, \alphat,\deltat)$. Indeed, if $|x_B|\le R$, from Lemma  $\ref{lema 3.5 de Gladis}$ we obtain that

\begin{equation*}
\|(1/w)\chi_B\|_{\infty}|B|^{\frac{\alphat-\deltat}{n}-\frac{1}{r}}\left(\frac{1}{|B|}\int_B v^{r'}(y)dy\right)^{1/{r'}}\menort R^{\alphat-\deltat-n/r-\theta}\menort 1,
\end{equation*}
and, if $|x_B|> R$
\begin{align*}
\|(1/w)\chi_B\|_{\infty}|B|^{\frac{\alphat-\deltat}{n}-\frac{1}{r}}\left(\frac{1}{|B|}\int_B v^{r'}(y)dy\right)^{1/{r'}}&\menort R^{\alphat-\deltat-n/r}|x_B|^{-\theta}\menort 1,
\end{align*}
which prove the assertion.

Finally we consider the situation in $(\ref{cuarto})$. Let $\deltat=\alphat-n$ and let us first consider $1<r\le\infty$. Let $w=1$ and $v\in L^{r'}(\zR^n)$. Then, it is not difficult to prove that the pair $(v,w)$ belongs to $\H(r, \alphat,\deltat)$. Let us now take $r=1$ and the pair of weights $(v,w)$ with $w(x)=|x|^{\theta}+1$, with $\theta>0$ and $v(x)=e^{-|x|}$. Then,  simple computations show that $(v,w)\in \H(r, \alphat,\deltat)$.

We are done.
\end{dem}

\section{Proof of the main results}\label{sprimeraparte}
We now give some previous lemmas that we shall  use in the proofs of the main results.
\begin{lema}\label{acotacion de la I3}
Let $0\le \alpha <n$, $0<\delta<\min(\gamma, (n-\alpha)/m)$, $m\in \N\cup \{0\}$, and $1\le r\le \infty$. If $b\in \Lambda(\delta)$ and  $(v,w)\in \H(r,\alphat,\deltat)$, then
$$
\int_{(2B)^c} |b(x)-b(z)|^m |K_\alpha(x-z)-K_\alpha(y-z)||f(z)| dz \,\menort\, \|b\|^m_{\Lambda(\delta)}|B|^{\deltat/n}\left\|f/v\right\|_{r} \inf_Bw
$$
for all $x,y \in B$.
\end{lema}

\begin{dem}
Let us suppose that $1<r\le\infty$. The proof for the case $r=1$ follows similar arguments with the obvious changes.

Let $x,y \in B$. From the fact that  $b\in\Lambda(\delta)$ and $K_\alpha\in H^{*}_{\alpha, \infty}$, we have that
\begin{align*}
\int_{(2B)^c} & |b(x)-b(z)|^m |K_\alpha(x-z)-K_\alpha(y-z)||f(z)| dz\\
&\menort \|b\|^m_{\Lambda(\delta)}\int_{(2B)^c} |x-z|^{\delta m} |K_\alpha(x-z)-K_\alpha(y-z)||f(z)| dz\\
&\menort \|b\|^m_{\Lambda(\delta)}\sum_{j=1}^{\infty}\frac{2^{j\delta m}|B|^{\delta m/n+\gamma/n}}{2^{j(n-\alpha+\gamma)}|B|^{(n-\alpha+\gamma)/n}}\int_{2^{j+1}B\setminus 2^{j}B }|f(z)|dz\\
&\menort \|b\|^m_{\Lambda(\delta)}|B|^{\frac{\delta m-n+\alpha}{n}}\sum_{j=1}^{\infty}(2^{j})^{\delta m -n+\alpha-\gamma}\int_{2^{j+1}B\setminus 2^{j}B }|f(z)|vv^{-1}dz.
\end{align*}
We can now apply  H\"older's inequality to get
\begin{align*}
\int_{(2B)^c} &|b(x)-b(z)|^m |K_\alpha(x-z)-K_\alpha(y-z)||f(z)| dz\\
&\menort \|b\|^m_{\Lambda(\delta)}\left\|f/v\right\|_{r}|B|^{\frac{\delta m-n+\alpha}{n}}\sum_{j=1}^{\infty}(2^j)^{\delta m-n+\alpha-\gamma}\left(\int_{2^{j+1}B\setminus 2^{j}B }v^{r'}(z)dz\right)^{1/{r'}}\\
&\menort \|b\|^m_{\Lambda(\delta)}\left\|f/v\right\|_{r}|B|^{\frac{\delta}{n}}\sum_{j=1}^{\infty}(2^j)^{\delta-\gamma}\left(\int_{2^{j+1}B\setminus 2^{j}B }\frac{v^{r'}(z)}{|x_B-z|^{r'(n-\alphat+\delta)}}dz\right)^{1/{r'}}\\
&\menort \|b\|^m_{\Lambda(\delta)}\left\|f/v\right\|_{r}|B|^{\frac{\delta}{n}}\sum_{j=1}^{\infty}(2^j)^{\delta-\gamma}\left(\int_{\R^n\setminus B }\frac{v^{r'}(z)}{|x_B-z|^{r'(n-\alphat+\delta)}}dz\right)^{1/{r'}}
\end{align*}
Since $\sum_{j=1}^{\infty}(2^{j})^{\delta-\gamma}$ is finite, by using the global condition (\ref{parte global}), we get
\begin{align*}
\int_{(2B)^c} |b(x)-b(z)|^m |K_\alpha(x-z)-&K_\alpha(y-z)||f(z)| dz\\
&\menort \|b\|^m_{\Lambda(\delta)}|B|^{\frac{\deltat}{n}}\left\|f/v\right\|_{r}\inf_Bw.
\end{align*}

\end{dem}

It is easy to see that the following result is true.

\begin{lema}\label{acotacion de la I1}
Let $0\le \alpha <n$, $0<\delta<\min(\gamma, (n-\alpha)/m)$, $m\in \N\cup \{0\}$, and $1\le r\le \infty$. If $b\in \Lambda(\delta)$ and  $(v,w)\in \H(r,\alphat,\deltat)$, then
$$
\|(1/w)\chi_B\|_{\infty} \int_B |T^m_{\alpha,b}f\chi_{2B}(x)| dx \,\menort \|b\|^m_{\Lambda(\delta)}|B|^{1+\deltat/n}\left\|f/v\right\|_{L^r}
$$
for all ball $B  \in \R^n$.
\end{lema}

We now give the proofs of Theorems \ref{teo para integrales singulares} and \ref{teorema A}. We remark that the same proof holds for both theorems if we consider $0\le \alpha <n$.

\medskip

\begin{prueba}{Theorems \ref{teo para integrales singulares} and \ref{teorema A}}
Let $0\le \alpha <n$ and $f$
 such that $f/v \in L^r(\R^n)$. Let $B\subset\R^n$ be a ball and $x\in B$. We split  $f=f_1+f_2$ with $f_1=f\chi_{2B}$. Define $a_B= \frac{1}{|B|}\int_{B} T^{m}_{\alpha,b}f_2$. Then,
\begin{align*}
&\frac{\|(1/w)\chi_B\|_{\infty}}{|B|}\int_B|T^m_{\alpha,b}f(x)-a_B|\,dx\\ &\hspace{3cm}\menort \frac{\|(1/w)\chi_B\|_{\infty}}{|B|}\int_B|T^m_{\alpha,b}f_1(x)|\,dx + \frac{\|(1/w)\chi_B\|_{\infty}}{|B|}\int_B|T^m_{\alpha,b}f_2(x)-a_B|\,dx \\
&\hspace{3cm}= I_1+I_2
\end{align*}
By Lemma \ref{acotacion de la I1}, we have
\begin{equation}
I_1\menort \|b\|^m_{\Lambda(\delta)}|B|^{\deltat/n}\left\|f/v\right\|_{r}.
\end{equation}
For $I_2$, we first estimate the difference $|T^m_{\alpha,b}f_2(x)-a_B|$ for every $x\in B$. Since
$$
|T^m_{\alpha,b}f_2(x)-a_B|= |T^m_{\alpha,b}f_2(x)-(T^m_{\alpha,b}f_2)_B|\menort \frac{1}{|B|}\int_B |T^m_{\alpha,b}f_2(x)-T^m_{\alpha,b}f_2(y)|dy.
$$
Then,
\begin{equation}
I_2\menort \frac{\|(1/w)\chi_B\|_{\infty}}{|B|}\int_B\frac{1}{|B|}\int_B |T^m_{\alpha,b}f_2(x)-T^m_{\alpha,b}f_2(y)|dy\;dx.
\end{equation}
Let $A=|T^m_{\alpha,b}f_2(x)-T^m_{\alpha,b}f_2(y)|$. If $x,y \in B$
\begin{align*}
A &\menort \int_{(2B)^c}|(b(x)-b(z))^m K_\alpha (x-z)-(b(y)-b(z))^mK_\alpha (y-z)||f(z)|dz \\
&\menort \int_{(2B)^c}|b(x)-b(z)|^m |K_\alpha (x-z)-K_\alpha (y-z)||f(z)|dz\\
&+\int_{(2B)^c}|(b(x)-b(z))^m-(b(y)-b(z))^m ||K_\alpha (y-z)||f(z)|dz\\
&= I_3+I_4.
\end{align*}
By Lemma \ref{acotacion de la I3}, we have
$$
I_3\menort \|b\|^m_{\Lambda(\delta)}|B|^{\frac{\deltat}{n}}\left\|f/v\right\|_{r}\inf_Bw.
$$
In order to estimate $I_4$, we use that $b\in\Lambda(\delta)$ and $x,y \in B$. If $A_j=2^{j+1}B\setminus 2^jB $, then
\begin{align*}
I_4 &\menort |b(x)-b(y)|\sum_{k=0}^{m-1}\int_{(2B)^c}|b(x)-b(z)|^{m-1-k}|b(y)-b(z)|^k|K_\alpha(y-z)||f(z)|dz \\
& \menort \|b\|_{\Lambda(\delta)}^m |B|^{\delta/n}\sum_{j=1}^{\infty}|2^{j+1}B|^{{\delta(m-1)}/n}\int_{A_j}|K_\alpha(x-z)||f(z)|dz\\
& \menort \|b\|_{\Lambda(\delta)}^m |B|^{\delta/n}\sum_{j=1}^{\infty}\int_{A_j}\frac{|f(z)|}{|x_B-z|^{(n-\alphat+\delta)}}dz\\
& \menort \|b\|_{\Lambda(\delta)}^m |B|^{\delta/n}\int_{\R^n\setminus B}\frac{|f(z)|}{|x_B-z|^{(n-\alphat+\delta)}}dz.\\
\end{align*}
Then, by H\"older's inequality, and the global condition (\ref{parte global}), we have

\begin{align*}
I_4 & \menort \|b\|_{\Lambda(\delta)}^m |B|^{\delta/n}\|f/v\|_{r}\left(\int_{\R^n\setminus B}\frac{v^{r'}(z)}{|x_B-z|^{r'(n-\alphat+\delta)}}dz\right)^{1/{r'}}\\
& \menort \|b\|_{\Lambda(\delta)}^m |B|^{\deltat/n}\|f/v\|_{r}\inf_Bw.
\end{align*}
Thus we get
\begin{align*}
I_2 & =  \frac{\|(1/w)\chi_B\|_{\infty}}{|B|} \int_B \frac{1}{|B|}\int_B (I_3+I_4)dydx\\
&\menort \|b\|^m_{\Lambda(\delta)}|B|^{\frac{\deltat}{n}}\left\|f/v\right\|_{r}\\
\end{align*}
and then, we have
$$
 \|(1/w)\chi_B\|_{\infty}\int_B|T^m_{\alpha,b}f(x)-a_B|\,dx\menort \|b\|^m_{\Lambda(\delta)}|B|^{\frac{\deltat}{n}+1}\left\|f/v\right\|_{r},
$$
so it remains to take supremum over all the balls $B$ to get the desired result.
\end{prueba}

\section{Preliminaries for extrapolation results}\label{extra1}
In this section we give some preliminary results in order to prove
Theorem $\ref{teo2}$. The proof of the first lemma is
straightforward and we omit it.

\begin{lem}\label{P1}
	Let $p(\cdot), r(\cdot) \in \mathcal{P}({\zR}^{n})$. If ${u}\in A_{p(\cdot),\infty}$ then ${u}\in A_{p(\cdot),\, r(\cdot)}$.
\end{lem}

The following lemma establishes a Reverse H\"{o}lder property under
certain conditions on the exponents involved.

\begin{lem}[\cite{MP2}]\label{luchi-Pradolini}
	
	Let $p(\cdot), r(\cdot) \in \mathcal{P}^{log}(\zR^n)$ such that $r^+\leq p^-$. Suppose that ${1}/{r(\cdot)}={1}/{p(\cdot)}+{1}/{s(\cdot)}$. Then
	\begin{equation}\label{equiv}
	 \|\chi_{Q}\|_{r(\cdot)}\simeq{\|\chi_{Q}\|_{s(\cdot)}}{\|\chi_{Q}\|_{p(\cdot)}}
	\end{equation}
\end{lem}
{
	\begin{rem}
		When
		$s(\cdot)=p'(\cdot)$ and, consequently $r(\cdot)=1$, a proof can be
		found in \cite{DHHR}.
	\end{rem}
	\begin{rem}\label{equiv1}It is easy to see that inequality \eqref{equiv} can be also written as
		\begin{equation*}
		 \|\chi_{Q}\|_{p'(\cdot)}\simeq{\|\chi_{Q}\|_{s(\cdot)}}{\|\chi_{Q}\|_{r'(\cdot)}}.
		\end{equation*}
	\end{rem}

\begin{lem}\label{first}
	Let $p(\cdot), r(\cdot) \in \mathcal{P}^{log}(\zR^n)$ such that $r^+\leq p^-$. Suppose that ${1}/{r(\cdot)}={1}/{p(\cdot)}+{1}/{s(\cdot)}$. If ${u}\in A_{r(\cdot), \infty}$ then ${u}\in A_{p(\cdot), \infty}$.
\end{lem}

\begin{dem}
	Since $1/p'(\cdot)=1/r'(\cdot)+1/s(\cdot)$, from H\"{o}lder's inequality, the hypothesis on the weight and remark $\ref{equiv1}$ we obtain that
	\begin{eqnarray*}
		 \|{u}^{-1}\chi_B\|_{p'(\cdot)}\|{u}\chi_B\|_{\infty}&\le& C \|{u}^{-1}\chi_B\|_{r'(\cdot)}\|{u}\chi_B\|_{\infty}\|\chi_B\|_{s(\cdot)}\\
		&\le& \|\chi_B\|_{r'(\cdot)}\|\chi_B\|_{s(\cdot)}\\
		&\le& C\|\chi_B\|_{p'(\cdot)}.
	\end{eqnarray*}
\end{dem}

\begin{lem}\label{P5}
	Let $s>1$, $p(\cdot)\in \mathcal{P}^{log}(\zR^n)$ such that $1<p^-\le p^+<s$ and $1/r(\cdot)=1/p(\cdot) - 1/s$. If ${u}\in A_{p(\cdot),\,r(\cdot)}$ then ${u}^{-s'}\in A_{p'(\cdot)/s'}$.
\end{lem}
\begin{dem}Let $Q\subset \zR^{n}$ be a cube. Since the inequality
	\begin{equation*}
	 {\|{u}^{-s'}\chi_{Q}\|_{p'(\cdot)/s'}\|{u}^{s'}\chi_{Q}\|_{(p'(\cdot)/s')'}}\leq C{|Q|}
	\end{equation*}
	is equivalent to
	\begin{equation*}
	 {\|{u}^{-1}\chi_{Q}\|_{p'(\cdot)}\|{u}\chi_{Q}\|_{(p'(\cdot)/s')'s'}}\leq C{|Q|^{1/s'}},
	\end{equation*}
	we prove the last one.

	Since $r(\cdot)=(p'(\cdot)/s')'s'$ and $1/s'=1/r(\cdot) ~+~1/p'(\cdot) $, from Lemma $\ref{luchi-Pradolini}$ we get that $|Q|^{1/s'}=\|\chi_{Q}\|_{s'}\simeq \|\chi_{Q}\|_{r(\cdot)}\|\chi_{Q}\|_{p'(\cdot)}$. By the hypothesis on the weight ${u}$ we obtain that
	
	\begin{eqnarray*}
		 \frac{\|{u}^{-1}\chi_{Q}\|_{p'(\cdot)}\|{u}\chi_{Q}\|_{(p'(\cdot)/s')'s'}}{|Q|^{1/s'}}&\leq& C \frac{\|{u}\chi_{Q}\|_{r(\cdot)}\|{u}^{-1}\chi_{Q}\|_{p'(\cdot)}}{\|\chi_{Q}\|_{r(\cdot)}\|\chi_{Q}\|_{p'(\cdot)}}\\
		&\le&C.
	\end{eqnarray*}
	We are done.
\end{dem}

\begin{prop}\label{P6}
	Let $p(\cdot)\in \mathcal{P}^{log}(\zR^n)$ and $s$, ${ \nu}$  such that
	$1 <p^-\le p^+<s$ and $1<{ \nu} < \min\{p^{-}, 2\}$. Let $1/r(\cdot)=1/p(\cdot) - 1/s$. If
	${u^{ \nu}}\in A_{p(\cdot)/{ \nu}, \infty}$ then ${ u}\in A_{r(\cdot),\infty}$.
\end{prop}
\begin{dem}
	Let $Q\subset \zR^{n}$ be a cube. Since $1<{ \nu} < p_{-}\leq
	p(\cdot)\leq r(\cdot)$, it is easy to check that there exists
	$\alpha >0$ such that
	 ${\nu}/{r'(\cdot)}={1}/{(r(\cdot)/{\nu})'}+{1}/{\alpha}$. By the
	generalized H\"{o}lder inequality and Lemma \ref{luchi-Pradolini} we
	get
	\begin{align*}
	\left(\frac{\|\chi_{B}{u}\|_{\infty} \|\chi_{B}{u}^{-1}\|_{r'(\cdot)}}{\|\chi_{B}\|_{r'(\cdot)}}\right)^{\nu}
	&\leq \frac{\|\chi_{B}{u^{\nu}}\|_{\infty} \|\chi_{B}{u^{-\nu} }\|_{r'(\cdot)/{\nu}}}{\|\chi_{B}\|_{r'(\cdot)/{\nu}}}\\
	&\leq C \frac{\|\chi_{B}{u^{\nu}}\|_{\infty} \|\chi_{B}{u^{-\nu}}\|_{(r(\cdot)/{\nu})'}\|\chi_{B}\|_{\alpha}}{\|\chi_{B}\|_{r'(\cdot)/{\nu}}}\\
	&\leq C \frac{\|\chi_{B}{u^{\nu}}\|_{\infty}
		 \|\chi_{B}{u^{-\nu}}\|_{(r(\cdot)/{\nu})'}}{\|\chi_{B}\|_{(r(\cdot)/{\nu})'}}\\
	&\leq
	C,
	\end{align*}
	{where we have used Lemma \ref{first} and the hypothesis on the weight to conclude that ${u^{\nu}}\in A_{r(\cdot)/{\nu}, \infty}$ }
\end{dem}

The following result can be found in \cite{CUW}.

\begin{prop}[\cite{CUW}]\label{P2}
	Given $q(\cdot)\in \mathcal{P}({\zR}^{n}) $, supposse that $\mu $ is a weight such that the Hardy-Littlewood maximal operator $M$ is bounded on $L^{q(\cdot)}(\mu)$. Fixed constants
	$\alpha > 0$, $\gamma\in\zR$ and a weight $\nu$, let $v$ be the weight defined by $v=\nu^{\gamma / \alpha}\mu^{1/\alpha}$. If  ${h}\in L^{\alpha q(\cdot)}(v) $,  then there exists ${H}\in L^{\alpha q(\cdot)}(v)$ such that:
	
	\noindent $(1)$ ${H}\geq {h}$;
	
	\noindent $(2)$ $\|{H}v\|_{\alpha q(\cdot)}\leq 2 \|{h}v\|_{\alpha q(\cdot)}$,
	
	\noindent $(3)$ ${H}^{\alpha}\nu^{\gamma}\in A_{1}$, with $[{H}^{\alpha}\nu^{\gamma}]_{A_{1}}\leq 2\|M\|_{L^{q(\cdot)}(\mu)}$.
\end{prop}

\section{Applications}\label{extra2}


In this section we give some applications of Theorem $\ref{teo2}$. We first state some previous results proved in \cite{MP} and
\cite{DPR} in order to extrapolate them, by means of Theorem $\ref{teo2}$, to the variable context.

We shall  suppose that $m\in \mathbb N_0=\mathbb N\,
\cup\{0\}$, and use the convention that $\beta/0=\infty$ if
$\beta>0$.

\subsection{Singular integral operators}

We consider two different types of operators $T$ according to the regularity condition satisfied by the kernel $K$.

\subsubsection{Singular integral operators with Lipschitz regularity}

We consider the singular integral operator defined in $(\ref{integral singular})$ with kernel $K$ satisfying conditions $S_0^*$ and $K_{0,\infty}^{*}$.

As a consequence of Corollary $\ref{luchigla}$ and Theorem \ref{teo2} with ${\nu} =1$ and $m\in \zN$, we get the following result.
\begin{thm}
	Let $0<\delta<\min\{\gamma,n/m\} $ and $b\in\Lambda(\delta)$. Let
	$p(\cdot)\in\mathcal{P}^{\log}(\zR^{n})$ and $0\leq
	\tilde{\delta}(\cdot)\leq \delta$ such that
	${1}/{p(\cdot)}+{\tilde{\delta}(\cdot)}/{n}={m\delta}/{n}$. If ${w^{-1}}\in
	A_{p(\cdot),\infty}$, then there exists a positive constant $C$ such
	that
	\begin{equation*}
	|||{T_{b}^{m}f}|||_{\mathbb{L}_{w}(\tilde{\delta}(\cdot))}\leq C
	\norm{b}_{\Lambda(\delta)}^m\norm{f/w}_{p(\cdot)}.
	\end{equation*}

\end{thm}

\subsubsection{Singular integral operators with H\"{o}rmander regularity}

We shall also be working with kernels $K$ satisfying certain
H\"ormander type regularity. Concretely, if $\Phi$ is a Young
function, we say that $K\in H_{\Phi}$  if there exist $c\ge 1$ and
$C>0$ such that for every $y\in \mathbb R^n$ and $R>c|y|$
\begin{equation}\label{Hphi}
\sum\limits_{j=1}^\infty (2^j R)^{n} \|\left(K(\cdot-y)-K(\cdot)\right)\chi_{|\cdot|\sim 2^j R}\|_{\Phi, B(0,2^{j+1}R)}  \leq C,
\end{equation}
where $|\cdot|\sim s$ means the set $\{x\in \mathbb R^n: s<|x|\leq 2s\}$.

For example, when $\Phi(t)=t^q$, $1\leq q<\infty$, we denote this class by
$H_{q}$ and it can be written as
\begin{equation*}
\sum\limits_{j=1}^\infty (2^j R)^{n}  \left(\frac{1}{(2^j
	R)^n}\int_{|x|\sim 2^j R}|K(x-y)-K(x)|^q
dx\right)^{1/q} \leq C.
\end{equation*}
We say that $K\in H_\infty$ if $K$ satisfies condition \eqref{Hphi}
with $\|.\|_{\Phi, B(0,2^{j+1}R)}$ replaced by
$\|.\|_{L^\infty,B(0,2^{j+1}R)}$.

The kernels given above are, a priori, less regular than the kernel
of the singular integral operator $T$ defined previously and they have been
studied by several authors. {For example, in \cite{LMRT}, the author
	studied singular integrals given by a multiplier. If
	$m:\mathbb{R}^n\rightarrow \mathbb{R}$ is a function, the multiplier
	operator $T_m$ is defined, through the Fourier transform, as
	$\widehat{T_mf}(\zeta)=m(\zeta) \widehat{f}(\zeta)$ for $f$ in the
	Schwartz class. Under certain conditions on the derivatives of $m$,
	the multiplier operator $T_m$ can be seen as the limit of
	convolution operators $T_m^N$, having a simpler form. Their
	corresponding kernels $K^N$ belong to the class $H_{r}$ with constant independent of $N$, for certain values
	of $r>1$ given by the regularity of the function $m$.}

The classes $H_q$, $1\leq q<\infty$, appeared implicit in \cite{KuW} where it is shown that the classical $L^q$-Dini condition for $K$ implies
$K\in H_q$ (see also \cite{RFRT} and \cite{WAT}).

Other examples of this type of operators are singular integrals operators
with rough kernels, that is, with kernel $K(x)=\Omega(x)|x|^{-n}$
where $\Omega$ is a function defined on the unit sphere $S^{n-1}$ of
$\mathbb{R}^n$, extended to $\mathbb{R}^n\setminus \{0\}$ radially.
The function $\Omega$ is an homogeneous function of degree $0$. In
\cite[Proposition~4.2]{LMRT}, the authors showed that $K\in
H_{\Phi}$, for certain Young function $\Phi$,
provided that $\Omega\in L^{\Phi}(S^{n-1})$ with
\[\int_0^1 \omega_\Phi(t) \frac{dt}{t}<\infty,\]
where $\omega_\Phi$ is the $L^\Phi$-modulus of continuity of $\Omega$ given by
\[\omega_\Phi(t)=\sup\limits_{|y|\leq t}\|\Omega(\cdot+y)-\Omega(\cdot)\|_{\Phi,S^{n-1}}<\infty,\]
for every $t\geq 0$.

Let $T^+$ be the differential transform operator studied in \cite{BLMRMT},
\cite{JoRo} and \cite{LMRT} and defined by
\begin{equation*}
T^{+}f(x)=\sum_{j\in \mathbb
	Z}(-1)^{j}\left(D_jf(x)-D_{j-1}f(x)\right),
\end{equation*}
where
\begin{equation*}
D_jf(x)=\frac{1}{2^j}\int_{x}^{x+2^j}f(t)\, dt
\end{equation*}
The operator above appears when dealing with the rate of convergence
of the averages $D_jf$, and it is a on-sided singular integral of
convolution type with a kernel $K$ supported in $(0, \infty)$ given
by
\begin{equation*}
K(x)=\sum_{j\in \mathbb
	 Z}(-1)^j\left(\frac{1}{2^j}\chi_{(-2^j,0)}(x)-\frac{1}{2^{j-1}}\chi_{(-2^{j-1},0)}(x)\right).
\end{equation*}
In \cite{LMRT} the authors proved that $K\in \bigcap_{r\ge 1}H_r$ but $K\notin H_\infty$. Moreover, $K\in H_{\psi}$, where $\psi(t)=\exp t^{1/1+\epsilon}-1$.

We shall be  dealing with conmmutators with symbols of
Lipschitz type. The smoothness condition associated to these
operators  is defined as follows.

\begin{defn}[\cite{DPR}]Let $m\in \mathbb N_0$, $0\leq \delta<\min\{1, n/m\}$ and let $\Phi$ be a Young function. We say that $K\in H_{\Phi,m}(\delta)$ if
	\begin{equation*}
	\sum\limits_{j=1}^\infty (2^j)^{m\delta}(2^j R)^{n} \|\left(K(\cdot-y)-K(\cdot)\right)\chi_{|\cdot|\sim 2^j R}\|_{\Phi, B(0,2^{j+1}R)}  \leq C.
	\end{equation*}
	for some constants $c\ge 1$ and $C>0$ and for every $y\in \mathbb   R^n$ with $R>c|y|$.
\end{defn}

Clearly, when $\delta=0$ or $m=0$, $H_{\Phi,m}(\delta)=H_{\Phi}$.

\begin{rem}\label{contentionsHdelta}
	It is easy to see that $H_{\Phi,m}(\delta_2)\subset H_{\Phi,m}(\delta_1) \subset H_{\Phi}$ whenever $0\leq \delta_1<\delta_2<\min\{1,n/m\}$.
\end{rem}

As we have mentioned above, Fourier multipliers and singular
integrals with rough kernels are examples of singular integral
operators with $K\in H_{\Phi}$ for certain Young function $\Phi$. By
assuming adequate conditions depending on $\delta$ on the multiplier
$m$, or on the $L^\Phi$-modulus of continuity $\omega_\Phi$, one can
obtain kernels $K\in H_{\Phi,m}(\delta)$, (see \cite{MP} for more
information).

{
	
	Related with the operators above, in \cite{MP} the authors proved
	the following result.

	\begin{thm}[\cite{MP}]\label{HormanderLipw}Let $0< \delta<\min\{ 1,n/m\}$ and $n/(m\delta)\leq s<n/((m-1)\delta)$. Let $w$ be a weight such that $A_{s/{\nu},\,\infty}$ for some $1<{\nu}<s$. Assume
		that $T$ has a kernel $K\in H_{\Phi,m}(\delta)$ for a Young function $\Phi$ such that $\Phi^{-1}(t)\lesssim t^{\frac{{\nu}-1}{s}}$ for every $t>0$. If
		$b\in\Lambda(\delta)$, then there exists a positive constant $C$
		such that
		\begin{equation*}
		|||{T_{b}^{m}f}|||_{\mathbb{L}_{w}(m\delta
			-{n}/{s})}\leq C \norm{b}_{\Lambda(\delta)}^{m}\norm{f{/w}}_{L^{s}}
		\end{equation*}
		for every $f\in L^s_w(\mathbb R^n)$.
	\end{thm}

	Then, by applying Theorem \ref{teo2} we obtain the following result.
	
	\begin{thm}
		Let $0< \delta<\min\{ 1,n/m\}$ and $n/(m\delta)\leq
		s<n/((m-1)\delta)$. Let $p(\cdot)\in \mathcal{P}^{\log}(\zR^{n})$
		and $0\leq \tilde{\delta}(\cdot)\leq \delta$ such that
		${m\delta}/{n}={1}/{p(\cdot)}+{\tilde{\delta}(\cdot)}/{n}$. Let
		$1<{\nu}<s$ and assume that $T$ has a kernel $K\in
		H_{\Phi,m}(\delta)$ for a Young function $\Phi$ such that
		$\Phi^{-1}(t)\lesssim t^{\frac{{\nu}-1}{s}}$ for every $t>0$. If
		$b\in\Lambda(\delta)$, then there exists a positive constant $C$
		such that
		\begin{equation*}
		\|T^{m}_{b}f\|_{\mathbb{L}_{w}(\tilde{\delta}(\cdot))}\leq C\|b\|_{\Lambda(\delta)}^{m}\|f/w\|_{p(\cdot)}
		\end{equation*}
		for any weight $w$ such that ${w^{-\nu}}\in A_{\frac{p(\cdot)}{\nu},\infty}$.
	\end{thm}

\subsection{Fractional type operators}

We now consider the commutators of fractional integral operator $T_{\alpha}$ and we classify them into two different
types, according to the smoothness conditions satisfied by $K_\alpha$.

\subsubsection{Fractional integral operators with Lipschitz regularity}

We first supposse that $K_{\alpha}\in H_{\alpha, \infty}^*$. As a consequence of Corollary $\ref{corolario fraccionario}$ and Theorem \ref{teo2}, with ${\nu} =1$, we obtain the following result.

\begin{thm}\label{aplicacion1}Let $0<\alpha<n$, and
	$0<\delta<\min\{\gamma,(n-\alpha)/m\} $. Let $n/(m\delta+\alpha)\le
	s<n/(\alpha+(m-1)\delta)$, if $m\in \mathbb{N}$ or $n/\alpha\le
	s<n/(\alpha-\gamma)_+$, if $m=0$. Let $p(\cdot)\in \mathcal{P}^{\log}(\zR^{n})$, ${{\tilde{\delta}}(\cdot)}/{n}= {(m\delta +\alpha)}/{n}-{1}/{p(\cdot)}$ with
	${0\le\tilde{\delta}(\cdot)}\le m\delta +\alpha-{n}/{s}$ and $b\in\Lambda(\delta)$. Then, if ${w^{-1}}\in A_{p(\cdot),\infty}$
	\begin{equation*}
	|||T^{m}_{\alpha,b}f|||_{\mathbb{L}_{w}({\tilde{\delta}}(\cdot))}\leq C\|b\|^{m}_{\Lambda (\delta)}\|f/w\|_{p(\cdot)}.
	\end{equation*}
	
\end{thm}

\subsubsection{Fractional integral operators with H\"ormander type
	regularity}\label{subsechormander}

We now introduce the smoothness condition on the kernel $K_\alpha$ that will be considered in this section.

We say that $K_{\alpha}\in H_{\alpha, \Phi}$ if
there exist $c\ge 1$ and $C>0$ such that for every $y\in \mathbb
R^n$ and $R>c|y|$
\begin{equation*}
\sum\limits_{j=1}^\infty (2^j R)^{n-\alpha} \|\left(K_{\alpha}(\cdot-y)-K_{\alpha}(\cdot)\right)\chi_{|\cdot|\sim 2^j R}\|_{\Phi, B(0,2^{j+1}R)}  \leq C,
\end{equation*}
where $|\cdot|\sim s$ means the set $\{x\in \mathbb R^n: s<|x|\leq 2s\}$.

When $\Phi(t)=t^q$, $1\leq q<\infty$, we denote this class by
$H_{\alpha,q}$ and it can be written as
\begin{equation*}
\sum\limits_{j=1}^\infty (2^j R)^{n-\alpha}  \left(\frac{1}{(2^j
	R)^n}\int_{|x|\sim 2^j R}|K_{\alpha}(x-y)-K_{\alpha}(x)|^q
dx\right)^{1/q} \leq C.
\end{equation*}

As in the case of the singular integral operators, the kernels given
above are less regular than the kernel of the fractional integral
operator $I_\alpha$. In \cite{Ku}, the author studied fractional
integrals given by a multiplier defined as in the case of singular
integral operators. Under certain conditions on the derivatives of
$m$, the multiplier operator $T_m$ can be seen as the limit of
convolution operators $T_m^N$, having a simpler form. Their
corresponding kernels $K_\alpha^N$ belong to the class $S_\alpha\cap
H_{\alpha,r}$ with constant independent of $N$, for certain values
of $r>1$ given by the regularity of the function $m$ (see
\cite{Ku}).

Other examples of this type of operators are fractional integrals
with rough kernels, that is, with kernel $K_\alpha(x)=\Omega(x)|x|^{\alpha-n}$
where $\Omega$ is a function defined on the unit sphere $S^{n-1}$ of
$\mathbb{R}^n$, extended to $\mathbb{R}^n\setminus \{0\}$ radially.
The function $\Omega$ is an homogeneous function of degree $0$. In
\cite[Proposition~4.2]{BLR}, the authors showed that $K_\alpha\in
S_\alpha\cap H_{\alpha,\Phi}$, for certain Young function $\Phi$,
provided that $\Omega\in L^{\Phi}(S^{n-1})$ with
\[\int_0^1 \omega_\Phi(t) \frac{dt}{t}<\infty,\]
where $\omega_\Phi$ is the $L^\Phi$-modulus of continuity given by
\[\omega_\Phi(t)=\sup\limits_{|y|\leq t}||\Omega(\cdot+y)-\Omega(\cdot)||_{\Phi,S^{n-1}}<\infty,\]
for every $t\geq 0$. This type of operators where also studied in \cite{CWW} and \cite{DL}.

Since we are dealing with symbols of Lipschitz type, the smoothness condition associated to the commutators  of $T_{\alpha}$ is defined as follows.

\begin{defn}Let $m\in \mathbb N_0$, $0<\alpha<n$, $0\leq \delta<\min\{1, (n-\alpha)/m\}$ and let $\Phi$ be a Young function. We say that $K_\alpha\in H_{\alpha,\Phi,m}(\delta)$ if
	\begin{equation*}
	\sum\limits_{j=1}^\infty (2^j)^{m\delta}(2^j R)^{n-\alpha} \|\left(K_{\alpha}(\cdot-y)-K_{\alpha}(\cdot)\right)\chi_{|\cdot|\sim 2^j R}\|_{\Phi, B(0,2^{j+1}R)}  \leq C.
	\end{equation*}
	for some constants $c\ge 1$ and $C>0$ and for every $y\in \mathbb   R^n$ with $R>c|y|$.
\end{defn}

Clearly, when $\delta=0$ or $m=0$, $H_{\alpha,\Phi,m}(\delta)=H_{\alpha,\Phi}$.

Recall that Fourier multipliers and fractional integrals with rough kernels are examples of fractional integral operators with $K_\alpha\in H_{\alpha,\Phi}$ for certain Young function. By assuming adequate conditions depending on $\delta$ on the multiplier $m$, or on the $L^\Phi$-modulus of continuity $\omega_\Phi$, we can obtain kernels $K_\alpha\in H_{\alpha,\Phi,m}(\delta)$. This fact can be proved by adapting Proposition 4.2 and Corollary 4.3 given in \cite{BLR} (see also \cite{LMRT}).

The following result is a generalization of Corollary \ref{corolario fraccionario} proved in \cite{DPR}. We shall consider again $m\in \mathbb N\cup \{0\}$.

\begin{thm}[\cite{DPR}]\label{HormanderLipw1}Let $0<\alpha <n$,
	$0< \delta<\min\{ 1,(n -\alpha)/m\}$ and $n/(m\delta+\alpha)\leq s<n/((m-1)\delta+\alpha)$. Let $w$ be a weight such that ${ w^{-\nu}}\in A_{s/{\nu},\infty}$ for some $1<{\nu}<s$. Assume
	that $T_\alpha$ has a kernel $K_\alpha\in S_\alpha\cap H_{\alpha,
		\Phi,m}(\delta)$ for a Young function $\Phi$ such that $\Phi^{-1}(t)\lesssim t^{\frac{{\nu}-1}{s}}$ for every $t>0$. If
	$b\in\Lambda(\delta)$, then there exists a positive constant $C$
	such that
	\begin{equation*}
	|||T_{\alpha,b}^{m}f|||_{\mathbb{L}_{w}(m\delta
		+\alpha-{n}/{s})}\leq C \norm{b}_{\Lambda(\delta)}^{m}\norm{f/w}_{L^{s}}
	\end{equation*}
	for every $f\in L^s_w(\mathbb R^n)$.
\end{thm}
Then, from the theorem above and as an application of Theorem $\ref{teo2}$ we obtain the following result. We shall now consider $m\in \zN $.

\begin{thm}\label{aplicacion3}
	Let $0<\alpha <n$, $0< \delta<\min\{ 1,(n -\alpha)/m\}$ and\, $n/(m\delta+\alpha)\leq s<n/((m-1)\delta+\alpha)$.
	Let $p(\cdot)\in \mathcal{P}^{\log}(\zR^{n})$  and ${{\tilde{\delta}}(\cdot)}/{n}={1}/{\beta}-{1}/{p(\cdot)}$, where $\beta = {n}/{(m\delta +\alpha)}$ and $0\le {\tilde{\delta}}(\cdot)\le m\delta+\alpha-n/s$. Let $1<{\nu}<s$ and assume
	that $T_\alpha$ has a kernel $K_\alpha\in S_\alpha\cap H_{\alpha,
		\Phi,m}(\delta)$ for any Young function $\Phi$ such that $\Phi^{-1}(t)\lesssim t^{\frac{{\nu}-1}{s}}$ for every $t>0$. If
	$b\in\Lambda(\delta)$ then, for every weight $w$ such that ${w^{-\nu}}\in A_{{p(\cdot)}/{{\nu}},\infty}$, there exists a positive constant $C$ such that
	\begin{equation*}
	|||T^{m}_{\alpha,b}f|||_{\mathbb{L}_{w}({\tilde{\delta}}(\cdot))}\leq C\|b\|^{m}_{\Lambda (\delta)}\|f\|_{L^{p(\cdot)}(w)},
	\end{equation*}
\end{thm}

\medskip

In order to give the next  result proved in \cite{DPR}, we introduce
some previous notation.

Given $0<\alpha<n$, $1\leq \beta<p<n/\alpha$ and a Young function $\Psi$, we
shall say that $\Psi\in \mathcal B_{\alpha,\beta}$ if
$t^{-\alpha/n}\Psi^{-1}(t)$ is the inverse of a Young function and
$\Psi^{1+{\rho\alpha}/{n}}\in B_{\rho}$ for every
$\rho>n\beta/(n-\alpha\beta)$, that is, there exists a positive constant $c$
such that
\[\int_c^\infty \frac{\Psi^{1+\frac{\rho\alpha}{n}}(t)}{t^\rho} \frac{dt}{t}<\infty\]
for each of those values of $\rho$.

For $k=0,1,...,m$, $m\in \mathbb N$, we denote
$c_{k}={m!}/{(k!(m-k)!)}$. If also $x$, $u\in \mathbb R^n$, we denote $S
(x,u,k)=(b(x)-b_{B})^{m-k}T_{\alpha}((b-b_{B})^{k}f_{2})(u)$, where $f_2=f\chi_{\mathbb R^n\setminus B}$ for a given ball $B$ and a locally integrable function $f$.

\begin{thm}[\cite{DPR}]\label{limitcaseH}Let $0< \delta<\min\{1,(n-\alpha)/m\}$ and let $s$ define by  $s= n/((m-1)\delta+\alpha)$. Let $w$ be a weight such that ${w^{-\nu}} \in A_{s/{\nu},\infty}$ for some $1<{\nu}<s$. Let $T_\alpha$ be a fractional integral operator with kernel $K_\alpha\in S_{\alpha} \cap H_{\alpha,\Phi,m}(\delta)$ where $\Phi$ is a Young function verifying $\Phi^{-1}(t)\lesssim t^{\frac{{\nu}-1}{s}}$ for every $t>0$, and $\widetilde{\Phi}\in \mathcal B_{m\delta+\alpha,{\nu}}$. If $b\in \Lambda(\delta)$,
	the following statements are equivalent,
	\begin{enumerate}
		\item \label{ecu1H} $|||T^{m}_{\alpha,b}f|||_{\mathbb{L}_{w}({\delta})}
	\leq C \|f{/w}\|_{s}$;
		\item There exists a positive constant $C$ such that
		\begin{equation}\label{ecu2H2}
		 \frac{\|{w^{-1}}\chi_B\|_\infty}{|B|^{1+\frac{\delta}{n}}}
		 \int_B\abs{\sum_{k=0}^{m}c_{k}\left[S(x,u,k)-(S(\cdot,u,k))_{B}\right]}dx\leq C \norm{f{/w}}_{s},
		\end{equation}
		for every ball $B\subset \mathbb{R}^{n}$, $x, u\in B$ and $f\in L^{1}_{{\rm loc}}(\mathbb R^n)$.
	\end{enumerate}
\end{thm}

Then, by applying Theorem \ref{teo2} we obtain the following result.

\begin{thm}\label{alpicacion4}
	Let $0< \delta<\min\{1,(n-\alpha)/m\}$ and $s=
	n/((m-1)\delta+\alpha)$. Let $1<{\nu}<s$ and assume that $T_\alpha$
	is a fractional integral operator with kernel $K_\alpha\in
	S_{\alpha} \cap H_{\alpha,\Phi,m}(\delta)$ where $\Phi$ is a Young
	function verifying $\Phi^{-1}(t)\lesssim t^{\frac{{\nu}-1}{s}}$ for
	every $t>0$, and $\widetilde{\Phi}\in \mathcal
	B_{m\delta+\alpha,{\nu}}$. Let
	$p(\cdot)\in\mathcal{P}^{\log}(\zR^{n})$ and $0\leq \tilde{\delta
	}(\cdot)\leq \delta <1$ such that
	 ${1}/{p(\cdot)}+{\tilde{\delta}(\cdot)}/{n}={1}/{\beta}={1}/{s}+{\delta}/{n}={m\delta
		+\alpha}/{n}$. Let $b\in \Lambda(\delta)$. If \eqref{ecu2H2} holds
	then
	\begin{equation*}
	|||T^{m}_{\alpha,b}f|||_{\mathbb{L}_{w}(\tilde{\delta}(\cdot))}
	\leq C \|f{w^{-1}}\|_{p(\cdot)},
	\end{equation*}
	for ${w^{-\nu}}\in A_{{p(\cdot)}/{{\nu}},\infty}$.
\end{thm}

\section{Proof of the main extrapolation result  }\label{extra3}

\begin{prueba}{Proposition \ref{P4}}
	Since ${\delta(\cdot)}/{n}={1}/{\beta}-{1}/{s(\cdot)}$, from the
	definition of $f^{\sharp}_{\delta(\cdot)}$ we have that
	\begin{equation*}
	 \frac{1}{|B|^{1/\beta}}\frac{1}{\|\chi_{B}\|_{s'(\cdot)}}\int_{B}|f(x)-m_{B}f|\leq
	f^{\sharp}_{\delta(\cdot)}(z),
	\end{equation*}
	for a.\,e. $z\in B$. By integrating over $B$ and applying
	H\"{o}lder's inequality we obtain that
	\begin{eqnarray*}
		 \frac{|B|}{|B|^{1/\beta}\|\chi_{B}\|_{s'(\cdot)}}\int_{B}|f(x)-m_{B}f|&\leq& \int_{B}f^{\sharp}_{\delta(\cdot)}(z)dz\\
		&\leq& C\|f^{\sharp}_{\delta(\cdot)}{w^{-1}}\|_{\frac{n}{\delta(\cdot)-\tilde{\delta}(\cdot)}}\|{w}\chi_{B}\|_{\frac{n}{n-(\delta(\cdot)-\tilde{\delta}(\cdot))}},
	\end{eqnarray*}
	and thus
	\begin{eqnarray*}
		 &&\frac{\|{w^{-1}}\chi_{B}\|_{\infty}}{|B|^{1/\beta}\|\chi_{B}\|_{p'(\cdot)}}
		\int_{B}|f(x)-m_{B}f| \\
		&&\hspace{3cm}\leq
		 \|f^{\sharp}_{\delta(\cdot)}{w^{-1}}\|_{\frac{n}{\delta(\cdot)-\tilde{\delta}(\cdot)}}
		 \|{w}\chi_{B}\|_{\frac{n}{n-(\delta(\cdot)-\tilde{\delta}(\cdot))}}\frac{\|\chi_{B}\|_{s'(\cdot)}\|{w^{-1}}\chi_{B}\|_{\infty}}{|B|~\|\chi_{B}\|_{p'(\cdot)}}.
	\end{eqnarray*}
	Since
	 ${n}/{(n-(\delta(\cdot)-\tilde{\delta}(\cdot)))}={1}/{s(\cdot)}+{1}/{p'(\cdot)}$,
	from the hypothesis on the weight and H\"{o}lder's inequality it
	follows that
	 $$\|{w}\chi_{B}\|_{\frac{n}{n-(\delta(\cdot)-\tilde{\delta}(\cdot))}}
	 \frac{\|\chi_{B}\|_{s'(\cdot)}\|{w^{-1}}\chi_{B}\|_{\infty}}{|B|~\|\chi_{B}\|_{p'(\cdot)}}\leq
	C.$$ Then, by taking supremun over $B$, we obtain the desired
	result.
\end{prueba}



\begin{prueba}{Theorem \ref{teo2}}
	Let ${1}/{r(\cdot)}={1}/{p(\cdot)}-{1}/{s}$, $\bar{r}={r}/{{\nu}}$, $\bar{p}={p}/{{\nu}}$ and $\bar{s}={s}/{{\nu}}$. Let $h$ be the function function defined by
	$$h=\left(\frac{|g|^{{\nu}} {w}^{{-\nu}\bar{p}'}}{\|g{w^{-1}}\|_{p(\cdot)}^{{\nu}}}\right)^{{\bar{p}}-{\bar{s}}}$$
	and $\tilde{h}=h^{-{\bar{s}\, '}/{\bar{s}}}$. Then $h\in L^{({p}/{s})'}({w}^{{{{\nu}\bar{p}'}}/{(p/s)'}})$ and $\|h\|_{L^{({p}/{s})'}({w}^{{{\nu}\bar{p}'}/{(p/s)'}})}\leq  1$. Indeed, since $\bar{p}'(p-{\nu})=p$
	\begin{align}\label{c1}
	 \nonumber\int_{\zR^{n}}h(x)^{(p(x)/s)\,'}{w(x)}^{{{\nu}\bar{p}'}}dx
	&= \int_{\zR^{n}}\left(\frac{|g(x)|{w(x)}^{{-}\bar{p}'(x)}}{\|g{w^{-1}}\|_{p(\cdot)}}\right)^{p(x)}{w(x)}^{{\nu}\bar{p}'}\\
	&= \int_{\zR^{n}}\frac{|g(x){w^{-1}(x)}|^{p(x)}}{\|g{w^{-1}}\|_{p(\cdot)}^{p(x)}}=
	1.
	\end{align}

	Moreover
	\begin{equation}\label{c4}
	\left(\int |g^{\nu}{w}^{{-\nu}\bar{p}'}|^{\bar{s}}h {w}^{{\nu} \bar{p}'}\right)^{1/\bar{s}}\leq \|g{w^{-1}}\|^{{\nu}}_{p(\cdot)},
	\end{equation}
	and, from ($\ref{c1}$),
	\begin{align*}
	 \int_{\zR^{n}}\tilde{h}(x)^{{\bar{r}}/{\bar{s}'}}{w(x)}^{{\nu}{\bar{p}'}}dx&=
	 \int_{\zR^{n}}h(x)^{-r(x)/s}{w(x)}^{{\nu}{\bar{p}'}}dx\\
	 &=\int_{\zR^{n}}h(x)^{(p(x)/s)'}{w(x)}^{{\nu}{\bar{p}'}}dx =1,
	\end{align*}
	that is $\tilde{h}\in L^{{\bar{r}}/{\bar{s}'}}({w}^{{\nu}{\bar{p}'\bar{s}'}/{\bar{r}}})$, $\|\tilde{h}\|_{L^{{\bar{r}}/{\bar{s}'}}({w}^{{\nu}{\bar{p}'\bar{s}'}/{\bar{r}}})}=1$.
	Now, from ($\ref{c4}$)
	\begin{equation}\label{a2}
	\|g{w^{-1}}\|^{{\nu}}_{p(\cdot)}\geq \left(\int
	 |g^{{\nu}}{w}^{{-\nu}\bar{p}'}|^{\bar{s}}\,\tilde{h}^{-{\bar{s}}/{\bar{s}'}}
	{w}^{{\nu} \bar{p}'}\right)^{1/\bar{s}},
	\end{equation}
	
	Note that, since ${1}/{\bar{r}}={1}/{\bar{p}}-{1}/{\bar{s}}$ and
	${w^{-\nu}}\in A_{\bar{p},\infty}$ then ${w^{-\nu}}\in A_{\bar{p},
		\bar{r}}$, and by Lemma $\ref{P5}$, ${w}^{{\nu} \bar{s}'}\in
	A_{{\bar{p}'}/{\bar{s}'}}$. Thus, the weight $\mu = {w}^{{-\nu}
		\bar{s}'}\in A_{({\bar{p}'}/{\bar{s}'})'}=A_{{\bar{r}}/{\bar{s}'}}$.
	If we now take $\alpha=1$ and
	$v={w}^{{\nu}{\bar{p}'\bar{s}'}/{\bar{r}}}$ in Proposition $\ref{P2}$,
	which implies that $\nu^{\gamma}={w}^{{\nu} \bar{p}'}$, then there
	exists $\tilde{H}\geq \tilde{h}$, such that
	$\|\tilde{H}v\|_{{\bar{r}}/{\bar{s}'}}\approx C$, which is
	equivalent to $\|\tilde{H}^{1/({\nu}
		\bar{s}')}{w}^{{\bar{p}'}/{\bar{r}}}\|_{r(\cdot)}\approx C$. Then,
	from \eqref{a2} we obtain that
	\begin{equation}\label{a5}
	\|g{w^{-1}}\|^{{\nu}}_{p(\cdot)}\geq \left(\int g^{s}
	\left(\tilde{H}^{-{1}/{({\nu}
			 \bar{s}')}}{w}^{{-}{\bar{p}'}/{\bar{s}'}}\right)^{s}\right)^{1/\bar{s}}
	\end{equation}
	
	On the other hand, from Proposition $\ref{P2}$ we also get that the
	weight $\tilde{H}\nu^{\gamma}=\tilde{H}{w}^{{\nu} \bar{p}'}\in A_{1}$,
	which is equivalent to $\tilde{H}^{-{1}/{ \bar{s}'}}{w}^{{{-\nu}
			\bar{p}'}/{\bar{s}'}}\in A_{\bar{s},\infty}$, that is
	$\left(\tilde{H}^{-\frac{1}{ {\nu}\bar{s}'}}{w}^{{-}\frac{
			\bar{p}'}{\bar{s}'}}\right)^{{\nu}}\in A_{\bar{s},\infty}$. Thus, by
	the hypothesis and from \eqref{a5} we obtain that
	\begin{align*}
	\|g{w^{-1}}\|_{p(\cdot)}&\geq |||f|||_{\mathbb{L}_{\tilde{H}^{-{1}/{({\nu} \bar{s}')}}{w}^{{-}{\bar{p}'}/{\bar{s}'}}}}(\delta) \\
	&\geq C \|f^{\sharp}_{\delta}\tilde{H}^{-{1}/{({\nu} \bar{s}')}}{w}^{{-}{\bar{p}'}/{\bar{s}'}}\|_{\infty}
	\|\tilde{H}^{{1}/{({\nu} \bar{s}')}}{w}^{{\bar{p}'}/{\bar{r}}}\|_{r(\cdot)}\\
	&\geq C \|f^{\sharp}_{\delta} {w}^{{-}\bar{p}'({1}/{\bar{s}'}-{1}/{\bar
			{r}})}\|_{r(\cdot)}\\
	&\geq C \|f^{\sharp}_{\delta} {w^{-1}}\|_{r(\cdot)}.
	\end{align*}
	From the hypothesis on $w$ and Proposition $\ref{P6}$ we conclude that ${w^{-1}}\in A_{r(\cdot), \infty}$,
	where $r(\cdot)=n/(\delta-\tilde{\delta}(\cdot))$. Thus, from Proposition $\ref{P4}$ we obtain that
	\begin{equation*}
	\|g{w^{-1}}\|_{p(\cdot)}\geq \|f^{\sharp}_{\tilde{\delta}(\cdot)}{w^{-1}}\|_{\infty}\cong |||f|||_{\mathbb{L}_{\tilde{\delta}(\cdot)}({w^{-1}})}.
	\end{equation*}
	
\end{prueba}

\small
\markright{}

\bibliographystyle{abbrv}

\noindent Gladis Pradolini, {\sl CONICET and Departamento de Matem\'atica, Facultad de Ingenier\'ia Qu\'imica, UNL, 3000, Santa Fe, Argentina}.\\
\noindent e-mail address: gladis.pradolini@gmail.com

\noindent Wilfredo A. Ramos  {\sl CONICET and Departamento de Matem\'atica, Facultad de Ciencias Exactas y Naturales y Agrimensura, UNNE, 3400, Corrientes, Argentina.}\\
\noindent e-mail address: oderfliw769@gmail.com

\noindent Jorgelina Recchi, {\sl CONICET and  Departamento de Matem\'aticas, UNS, 8000, Bah\'ia Blanca, Argentina.}\\
\noindent e-mail address: drecchi@uns.edu.ar, jrecchi@gmail.com
\end{document}